\documentclass{article}
\usepackage{eufrak}
\usepackage{euscript}
\usepackage{epsfig}
\usepackage{pictex}
\usepackage{graphicx}

\def\lnorm{\,\rule[-1mm]{0.6mm}{4mm}\,}

\def\1{\hbox{\upshape1\kern-.15em\vrule height 1.6ex width .3pt
\vrule width .8pt height .25pt\kern.15em}}
\def\der{\,\rule[0mm]{0.1mm}{1.3mm} \rule[0mm]{1.3mm}{0.1mm}
 \hspace{-1.3mm} \rule[1.3mm]{1.3mm}{0.1mm}
 \hspace{-1.3mm} \rule[2.3mm]{1.3mm}{0.1mm}
 \rule[0mm]{0.1mm}{2.3mm} \,}
\def\lev{\langle \,}
\def\des{\,\rangle}
\def\vol{{\rm vol}}

\newtheorem{lemma}{{\sc LEMMA}}[section]
\newtheorem{theorem}{{\sc THEOREM}}[section]
\newtheorem{proposition}{{\sc PROPOSITION}}[section]
\newtheorem{assumption}{Assumption}[section]
\newtheorem{procedure}{Discretization procedure}[section]
\newtheorem{definition}{DEFINITION}[section]
\newtheorem{corollary}{{\sc COROLLARY}}[section]
\newtheorem{example}{{\sc Example}}[section]

\newcommand {\Proof}  {\mbox{\sc Proof: }}
\newcommand {\QED}{\hfill {\bf \textsf{QED}}}

\newcommand {\Mg} {\overline{M}}
\newcommand {\Md} {\underline{M}}
\newcommand {\mess} {{\rm meas}\,}
\newcommand {\supp} {{\rm supp}\,}

\def\msb{\mbox{\boldmath{$b$}}}
\def\msx{\mbox{\boldmath{$x$}}}
\def\msxd{\mbox{\boldmath{\scriptsize$x$}}}
\def\msy{\mbox{\boldmath{$y$}}}
\def\msv{\mbox{\boldmath{$v$}}}
\def\msvd{\mbox{\boldmath{\scriptsize$v$}}}
\def\msw{\mbox{\boldmath{$w$}}}
\def\mswd{\mbox{\boldmath{\scriptsize$w$}}}
\def\msz{\mbox{\boldmath{$z$}}}
\def\mszd{\mbox{\boldmath{\scriptsize$z$}}}
\def\msp{\mbox{\boldmath{$p$}}}
\def\mspd{\mbox{\boldmath{\scriptsize$p$}}}
\def\msq{\mbox{\boldmath{$q$}}}
\def\msr{\mbox{\boldmath{$r$}}}
\def\msrd{\mbox{\boldmath{\scriptsize$r$}}}
\def\msR{\mbox{\boldmath{$R$}}}
\def\msRd{\mbox{\boldmath{\scriptsize$R$}}}
\def\msk{\mbox{\boldmath{$k$}}}
\def\mskd{\mbox{\boldmath{\scriptsize$k$}}}
\def\msm{\mbox{\boldmath{$m$}}}
\def\msmd{\mbox{\boldmath{\scriptsize$m$}}}

\def\msl{\mbox{\boldmath{$l$}}}
\def\msld{\mbox{\boldmath{\scriptsize$l$}}}
\def\mse{\mbox{\boldmath{$e$}}}
\def\msed{\mbox{\boldmath{\scriptsize$e$}}}
\def\msmu{\mbox{\boldmath{$\mu$}}}
\def\msch{\mbox{\boldmath{$\chi$}}}
\def\mst{\mbox{\boldmath{$t$}}}

\font\msbm = msbm10
\def\Msbm#1{\hbox{\msbm #1}}

\newcommand {\bbR}{{\Msbm R}}
\newcommand {\bbN}{{\Msbm N}}
\newcommand {\bbZ}{{\Msbm Z}}

\newcommand {\bbJ}{\1}

\newdimen\vskp
\def\AMSclass#1{\vskip\vskp\vbox{\hbox{\vbox{\def\ams{#1}%
	\ifx\ams\empty
	\errmessage{Please put AMS subject classification in \noexpand
	\AMSclass command}
	\else\noindent
	\small\strut{\bf AMS subject classification:\hskip0.5em}\ams
	\strut}}}\fi}

\def\KeyWords#1{\vskip\vskp\vbox{\hbox{\vbox{\def\Kwd{#1}%
	\ifx\Kwd\empty
	\errmessage{Please fill Key words in \noexpand
	\KeyWords command}
	\else
	\newdimen\keywdswidth
	\setbox0=\hbox{\small\bf Key words:\hskip0.5em}
	\keywdswidth=\wd0
	\noindent\hangindent=\keywdswidth
	\hangafter=1
	\small\strut{\bf Key words:\hskip0.5em}\Kwd\strut}}}\fi}

\title{Numerical approach to $L_1$-problems with the second order elliptic
operators\thanks{Supported
	by grant 037-1193086-2771 of the Ministry of Science, Higher Education and Sports,
Croatia.}}

\author{%
Ned\v zad Limi\'c
	\thanks{%
	Dept. of Mathematics, University of Zagreb, 
	Bijeni\v{c}ka 30, 10002 Zagreb, Croatia,
	e--mail: nlimic@math.hr} 
and Mladen Rogina	\thanks{%
	Dept. of Mathematics, University of Zagreb, 
	Bijeni\v{c}ka 30, 10002 Zagreb, Croatia,
	e--mail: rogina@math.hr} 
}
\begin{document}
\maketitle
\begin{abstract}
\noindent
For a second order differential operator $A(\msx) =-\nabla a(\msx)\nabla + 
b'( \msx)\nabla+ \nabla \big(\msb''( \msx) \,\cdot\big)$
on a bounded domain $D$ with the Dirichlet boundary conditions on
$\partial D$ there exists the inverse $T(\lambda, A)= (\lambda I+A)^{-1}$
in $L_1(D)$. If $\mu$ is a Radon (probability) measure
on Borel algebra of subsets of $D$, then $T(\lambda, A)\mu \in L_p(D),
p \in [1, d/(d-1))$. We construct the numerical approximations to 
$u =T(\lambda, A)\mu$ in two steps. In the first one we construct 
grid-solutions ${\bf u}_n$ and in the second step we embed
grid-solutions into the linear space of hat functions
$u(n) \in \dot{W}_p^1(D)$. The strong convergence to the original solutions
$u$ is established in $L_p(D)$ and the weak convergence in $\dot{W}_p^1(D)$.
\end{abstract}
\AMSclass{(2000) 35J20, 35J25, 35J15, 65N06, 65N15}

\section{Introduction}\label{sec1}
If the  tensor-valued function $\{a_{ij}\}_{11}^{dd}$ (diffusion tensor) on
a bounded domain $D \subset {\bbR}^d$ satisfies the strict ellipticity
conditions, then the second order differential operator in divergence form:
$A =  -\sum_{ij}\partial_i a_{ij}\partial_j + \sum_i b_i' \partial_i +
\sum_i \partial_i (b_i'' \cdot)$, with homogeneous Dirichlet boundary
condition on $\partial D$, has an inverse $A_D^{-1}$ in terms of an integral 
operator mapping the
Radon measures on Borelian sets $\mathfrak{B}(D)$ into $\dot{W}_p^1(D)$
for each $p \in [1, d/(d-1))$ \cite{BO,LR2}. There are partial results
on numerical solutions of the boundary value problem $A_Du = \mu$
\cite{Cl,LR2,LR3}. Here we extend results of \cite{LR3} to a general
boundary value problem on $D \subset {\bbR}^d$. The
operator $A_D$ is discretized by a {\em system matrix} $A_n$ and numerical
solutions are represented by grid-functions ${\bf u}_n$. Then ${\bf u}_n$
is imbedded into the linear space of hat functions and compared with
the solution of the original problem in order to prove convergence.
Our particular intention is to construct numerical approximations ${\bf u}_n$
with system matrices $A_n$ possessing {\em compartmental
structure}, that is a matrix having positive diagonal elements, non-positive 
off-diagonal elements and non-negative column sums. Apparently, such matrices
are transpose of M-matrices. Because
of our determination to look for approximations $A_n$ with the compartmental structure,
we have faced a number of non-typical problems in numerical analysis. All of
them stem from the construction of the numerical schemes, rather then from the structure of
convergence proofs.

Various definitions and auxiliary results which are necessary to prove
convergence are given in Sections~\ref{sec2} and \ref{sec3}.
In the case of dimension $d=2$ the proposed numerical scheme is 
already on a level of an algorithm and can be straightforwardly applied. The
construction of discretizations $A_n$ with the compartmental structure for
any dimension is carried out in Section \ref{sec4}. Grid-solutions of discretized
problems are imbedded in the linear space of hat functions and the
obtained approximate solutions are analyzed from the standpoint of convergence. 
The convergence in $W_2^1$-spaces is proved in Section \ref{sec5} and
the convergence in $W_p^1$-spaces in Section \ref{sec6}. Section~\ref{sec7} provides two
examples, demonstrating the efficiency of constructed schemes.
\section{Definition of the problem}\label{sec2}
Elements of ${\bbR}^d$ are denoted by $\msx, \msm,
\msp$ etc. The Euclidean norm in ${\bbR}^d$ is denoted by $|\cdot|$.
The only open subsets of ${\bbR}^d$ considered in this work,
are bounded and connected open sets with Lipshitz boundary \cite{Ma,Ste}. 
We call a subset of this kind a {\em domain with Lipshitz boundary}. 
We denote it by D, and its boundary
by $\partial D$. For a set $S \subset {\bbR}^d$ the closure is denoted by
$\overline{S}$ and sometimes by $cls(S)$.

The Banach spaces of functions $C^{(k)}({\bbR}^d), C^{(k)}(\overline{D})$ are defined
as usually. Their norms are denoted by $\Vert \cdot \Vert_\infty^{(k)}$. 
The closure of functions in $C^{(k)}({\bbR}^d)$ with compact supports
determines the subspace $C_0^{(k)}({\bbR}^d)$. The closure of
functions in $C(\overline{D})$ with supports in $D$ determines the
subspace $\dot{C}(\overline{D})$, and then $\dot{C}^{(k)}=C^{(k)}(\overline{D}) 
\cap \dot{C}(\overline{D})$.
The $L_p$-spaces as well as Sobolev $W_p^1$-spaces are defined in a
standard way \cite{Ma,Ste}.
Their norms are denoted by $\Vert \cdot \Vert_p$ and $\Vert \cdot \Vert_{p,1}$,
respectively. For each $p,\: 1 \leq p \leq \infty$, the norm of $W_p^1(D)$ is
defined by $\Vert u \Vert_{p,1} =  \big(\Vert u \Vert_p^2 + \Vert \nabla u 
\Vert_p^2\big)^{1/2}$, where $\Vert \,\nabla u\, \Vert_p \, = \, \big(\sum_{j=1}^d\,
\Vert \,\partial_j u\, \Vert_p^2\big)^{1/2}$. Because the domain $D$ has Lipshitz boundary,
the space $W_\infty^1(D)$ can be realized as the space of continuous
functions on $\overline{D}$, for which the first partial derivatives are
elements of $L_\infty(D)$. The completion in the norm of $W_p^1(D)$ of
functions in $C^{(1)}(D) \cap \dot{C}(\overline{D})$ is denoted by
$\dot{W}_p^1(D)$. The linear space $C(\overline{D}) \cap W_p^1(D)$ is dense in
$W_p^1(D)$ for all $p, \:1 \leq p \leq \infty$, and the closure of
$\dot{C}(\overline{D}) \cap W_p^1(D)$ in the norm of $W_p^1(D)$ is equal
to $\dot{W}_p^1(D)$ \cite{Ma}. Let $X$ be a Banach space and $X^\dag$ 
its dual. Then the value of $f \in X^\dag$ at $u \in X$ is denoted by
$\lev f | u\des$. Let ${\cal R}(D)$ 
be the convex set of positive Radon measures $\mu$ on $\mathfrak{B}(D)$.
Then $\langle v \,|\, \mu \rangle = \int_D v(\msx) \mu(d {\bf
x})$ is well defined for $v \in \dot{W}_\infty^1(D)$. We say that a sequence of
$\mu_n \in {\cal R}(D)$ converges weakly to $\mu \in {\cal R}(D)$ if $\lim_n \lev v |
\mu_n \des = \lev v | \mu \des$ for each $v \in \dot{C}(\overline{D})$.

Let ${\bbJ}_S$ be the {\em indicator} of $S \subset {\bbR}^d$, {\em i.e.}
${\bbJ}_S(\msx)=1$ for $\msx \in S$, and zero otherwise.
We say that $f$ on ${\bbR}^d$ is  {\em piecewise continuous with
respect to the decomposition} ${\bbR}^d =\cup_k D_k$ if there exist
a finite collection of $L$ disjoint, Lebesgue measurable subsets $D_k \subset {\bbR}^d$, and
bounded, uniformly continuous functions
on ${\bbR}^d$, $\{f_j\}_{j=1}^L$, such that
${\bbR}^d = \cup_{k=1}^L D_k$ and $f \ = \ \sum_{j = 1}^L f_j {\bbJ}_{D_j}$.
If not necessary a part of the terminology such as `with respect \ldots' is omitted.
A function $f_D$ on $D \subset {\bbR}^d$ is piecewise continuous if
there exists a piecewise continuous $f$ on ${\bbR}^d$ such that
$f_D = f |D$. Piecewise constant functions are special cases of piecewise
continuous ones.

We consider a $2^{{\rm nd}}$-order elliptic operator on ${\bbR}^d$,
\begin{equation}\label{exp2.1}
   A(\msx) = -\:\sum_{ i,j=1}^d  \partial_i a_{ij}(\msx) \partial_j
 \:+\: \sum_{j=1}^d b_j'(\msx)\, \partial_j
 \:+\: \sum_{j=1}^d \partial_j \big( b_j''(\msx) \cdot \big)\:+\: c(
 \msx),
\end{equation}
for which the coefficients must fulfill the following:
\begin{assumption}\label{Ass2.1}
The functions  $a_{ij}= a_{ji}$, $b_{i}', b_{i}''$ ($i,j = 1,2,\dots, d$) and $c$
are:
\begin{description}\itemsep 0.cm
 \item{a)} piecewise continuous on ${\bbR}^d$, $c \geq 0$ and $a_{ij}(
\msx)$ converge to constant values as $|\msx|$
increases,
 \item{b)} there are positive numbers $\underline{M}, \overline{M}, \:
0 < \underline{M} \leq \overline{M}$, such that the strict ellipticity 
\begin{equation}\label{exp2.2}
  \underline{M} \,|\msx|^2  \,\leq\,\sum_{i,j=1}^{d}  a_{ij}(\msx)
z_i \bar{z}_{j}
 \leq \,\overline{M}\, |\msx|^2 ,  \quad  \msx \in  {\bbR}^d
\end{equation}
holds.
\end{description}
\end{assumption}

The differential operator $A_0(\msx) = -\sum_{ i,j=1}^d
\partial_i a_{ij}(\msx) \partial_j$ is called
the {\em main part} of $A(\msx)$.

Let us define a real bilinear form on $W_q^1(D) \times W_p^1(D)$, $1/p+1/q = 1$, by:
\begin{equation}\label{exp2.4}\begin{array}{c}\displaystyle
 a(v,u) \ = \ \sum_{i,j =1}^d\: \int_D\: a_{ij}( \msx)\:
 \partial_i v(\msx)\:  \partial_j u(\msx) \:d\msx 
 - \sum_{i=1}^d \: \int_D\: b_i'(\msx) \:v(\msx)\,  \partial_i u(\msx) \:d \msx  
 \\ \displaystyle 
 - \sum_{i=1}^d \: \int_D\: b_i''(\msx) \:\partial_i v(\msx)\, 
 u(\msx) \:d \msx +  \int_D\: c(\msx)\, v(\msx)\,  u(\msx)\,d\msx .  \end{array}
\end{equation}
For a domain $D$ with Lipshitz boundary $\partial D$ and for each pair
$v \in \dot{W}_q^1(D)$, $u \in \dot{W}_p^1(D) \cap \{A u \,\in \, 
(\dot{W}_q^1(D))^{\dag}\} $, $1<p<\infty$, the Green formula must be valid,
\[ a(v,u) \ = \ \lev v\,|\,A\, u \des . \]
The Green formula is also valid for each pair
$v \in \dot{W}_{\infty}^1(D)$, $u \in \dot{W}_1^1(D)\cap \{ Au \in {\cal R}(D)\}$.

The boundary value problem, to be studied in this work, is defined by
\begin{equation}\label{exp2.3}\begin{array}{c}
 \big(\lambda I \,+\,A(\msx)\big)\, u(\msx) \: = \: \mu(\msx), \quad \msx \in D,\\
 u\,\big | \,\partial D \ = \ 0, \end{array}
\end{equation}
where $\lambda \geq 0$ and $D$ is a domain with Lipshitz boundary. In the case 
of $D = {\bbR}^d$ we suppose that $\lambda > 0$ and the boundary condition in 
(\ref{exp2.3}) is omitted. The nonhomogeneous term $\mu$ is a positive Radon measure 
for $p \in [1,d/(d-1))$, and $\mu \in W_2^{-1}(D)$ for $p = 2$.

The variational formulation of (\ref{exp2.3}) for a solution
$u \in \dot{W}_p^1(D)$, $p \in [1,d/(d-1))$ or $p=2$,  has the form:
\begin{equation}\label{exp2.9}
 \lambda (v\,|\,u)\,+\,a(v,u) \ = \ \langle \,v \,|\,\mu \,\rangle, \quad 
 \textrm{for~any} \ v \in \dot{W}_q^1(D).
\end{equation}
Solutions of (\ref{exp2.3}) and (\ref{exp2.9}) are called strong and weak
solutions, respectively. In the case of a problem on ${\bbR}^d$ the variational
problem is defined by expression
\begin{equation}\label{exp2.10}
 \lambda (v\,|\,u)\,+\,a(v,u) \ = \ \langle \,v \,|\,\mu \,\rangle, \quad 
 \textrm{for~any} \ v \in W_q^1({\bbR}^d),
\end{equation}
where for $p=2$ we have $\mu \in W_2^{-1}({\bbR}^d)$ and for $p \in [1,d/(d-1))$ 
$\mu \in {\cal R}(D)$ with a bounded $D \subset {\bbR}^d$.

For the differential operator $H(\msx) =\lambda I -\sigma^2\Delta$ on ${\bbR}^d$
and $\lambda > 0$ the fundamental solution $(\msx, \msy) \mapsto t(\lambda,\msx -\msy)$ 
can be represented in terms of the Bessel function $K_\nu, \nu = (d-2)/2$. Let us
denote the corresponding operator by $T(\lambda,H)$, i.e. we have
$(\lambda I -\sigma^2\Delta)T(\lambda,H) = I$ on the linear space of continuously
differentiable functions with compact supports. Much more, for any $\alpha > 0$ 
there exists a representation of the 
operator $T(\lambda,H)^\alpha$ in terms of an integral operator with a positive
kernel $t_\alpha(\lambda,\msx-\msy)$ which is expressed by the Bessel function 
$K_\nu,\nu = (d-2\alpha)/2$ \cite{Sh}. 
The Green function for the differential operator $\lambda I -\sigma^2\Delta$ on 
$D$ with the homogeneous Dirichlet boundary conditions at $\partial D$ is denoted 
by $t(\lambda,H_D,\cdot,\cdot)$ and the corresponding integral operator by 
$T(\lambda,H_D)$. In this way we have $(\lambda I -\sigma^2\Delta)T(\lambda,H_D)=
I$ on the linear space of continuously differentiable functions on $\overline{D}$.
There exist representations of $T(\lambda,H_D)^\alpha$ as integral operators
with kernels $t_\alpha(H_D,\cdot,\cdot)$ which are positive on $D \times D$.
The equality $(\lambda I+H_D)T(\lambda,D) = I$ enables us to define various closures
$H_D = -\lambda I+T(\lambda,H_D)^{-1}$, such as the closure from $\dot{W}_2^1(D)$
onto $W_2^{-1}(D)$, from $\mathfrak{D}(H) = T(\lambda,H_D)L_p(D)$ onto $L_p(D)$
{\em etc.}

For $\lambda$ sufficiently large and $\mu \in W_2^{-1}({\bbR}^d)$
solutions to~(\ref{exp2.3})-(\ref{exp2.10})
exist and can be represented by $T(\lambda,H)^{1/2}$ (or $T(\lambda,H_D
)^{1/2})$ as described in the following. Let us define the bounded operator
on $L_2({\bbR}^d)$ by
\begin{equation}\label{exp2.13}
 W \ = \ T(\lambda,H)^{1/2}\,\Big(\sum_{ij}\: \partial_i \,(a_{ij}
 \,-\,\Mg \delta_{ij}\,)\partial_j \Big)\,T(\lambda,H)^{1/2}.
\end{equation}
Then $\Vert W \Vert_2 \leq (1-\gamma)$ where $\gamma = \Md/\Mg < 1$, so that
there exist the operator
\begin{equation}\label{exp2.11}
 T(\lambda,A_0) \ = \ T(\lambda,H)^{1/2}\,(I-W)^{-1}\,T(\lambda,H)^{1/2}
\end{equation}
mapping $W_2^{-1}({\bbR}^d)$ into $W_2^1({\bbR}^d)$ with the norm
\begin{equation}\label{exp2.12}
 \Vert T(\lambda,A_0)\Vert_{\mathsf{L}(W_2^{-1},W_2^1)} \leq 
 \frac{\Mg}{\Md}\big(\lambda^{-1} \,+\,\Mg^{-1}\big).
\end{equation}
From Aronson's inequalities \cite{Ba} we have the following result.
The operator $T(\lambda,A_0)$ is an integral operator and its kernel is the
fundamental solution of differential operator $\lambda I+A_0(\msx)$ on 
${\bbR}^d$. If the lower order differential operators in (\ref{exp2.1}) 
are non-trivial, then $T(\lambda,A)$ exists for $\lambda$ sufficiently large. 
Results (\ref{exp2.11}) and (\ref{exp2.12}) are valid for bounded domains as well. 
We have to replace $H$ with the corresponding $H_D$ and obtain in this way
the operators $W_D, T(\lambda,A_D)$. Various closures $A, A_D$ are defined in
terms of the constructed operators $T(\lambda,A), T(\lambda,A_0)$, respectively,
as in the case of differential operator $-\sigma^2 \Delta$.

In the case of $\mu \in {\cal R}(D)$ solutions to (\ref{exp2.3})-(\ref{exp2.10})
also exist. We have the following result \cite{BO, LR2}:

\begin{theorem}\label{th2.1} Let $D$ be a bounded domain with Lipshitz
boundary. For each $p \in [1, d/(d-1))$ there exists
a unique weak solution $u$ of (\ref{exp2.9}) belonging to the class
$\dot{W}_p^1(D)$ and possessing the following properties:
\begin{description}\itemsep 0.cm
 \item{(i)} There exists a positive number $c$ depending on
$\underline{M}, \overline{M}, p, D$, such that the following inequality
is valid:
\[ \Vert \,u\, \Vert_{p,1} \ < \ c\,\mu(D).  \]
 \item{(ii)} If $\{ \mu_n : n \in {\bbN}\} \subset {\cal R}(D)$
converges weakly to a $\mu \in {\cal R}(D)$, then the corresponding sequence
of weak solutions $\{ u_n : n \in {\bbN}\} \subset \dot{W}_p^1(D)$, $u_n =
A_D^{-1}\mu_n$, converges strongly in $L_p(D)$ to $u = A_D^{-1} \mu$.
\end{description}
\end{theorem}

This theorem is the theoretical background for construction and analysis
of convergence of numerical solutions in $L_1({\bbR}^d)$.

\section{Grids and associated functions}\label{sec3}
Let the orthogonal coordinate system in ${\bbR}^d$ be determined by
unit vectors $\mse_i$, and let us define
the set $G_n$ by:
\begin{equation}\label{exp3.1}
 G_n \ = \ \{ \msx \:=\: h(n)\,\sum_{l=1}^d \:k_l\:
 \mse_l \ : \ k_l \in {\bbZ} \},
\end{equation}
where $h(n) = 2^{-n}$ is called the grid-step. A grid-step is usually denoted by
$h$ and only if necessary by $h(n)$. Elements of $G_n$ are
called grid-knots and the constructed sets $G_n, n \in {\bbN}$
are called grids. Sometimes we say that $G_n$ discretize 
${\bbR}^d$. Accordingly, the subgrids $G_n(D) \subset G_n$ defined by
$G_n(D) = G_n \cap D$ are called discretizations of $D$.
To each $\msv \in G_n$ there corresponds a grid-cube $C_n({\bf 1}, \msv) = 
\prod_1^d \,[v_j,v_j+h)$, where $v_j$ are coordinates of $\msv \in G_n$. 
Cubes $C_n({\bf 1},\msv)$ define a decomposition of 
${\bbR}^d$ into disjoint sets. Apart from the basic cubes, $C_n({\bf 1},\msv), 
\msv \in G_n$, for constructions we need larger sets. Let $\msp \in {\bbN}^d$. 
Then
\[ C_n(\msp,\msv) \ = \  \prod_{i=1}^d \: [v_i,v_i +hp_i)\]
are apparently rectangles with "lower left" vertices $\msv$
and edges of size $h p_i$. These rectangles define a partition of ${\bbR}^d$
as well. The considered cubes $C_n({\bf 1}, \msv)$ and
rectangles $C_n(\msp,\msv)$ are semi-closed
in the sense that they contain only one of their $2^d$ vertices. 

Basic cubes are defined by their "lower left" vertices. Apart from these
basic cubes for our constructions we need closed rectangles,
\begin{equation}\label{exp3.1a}
 S_n(\msp,\msv) \ = \ \prod_{i=1}^d \: [v_i - hp_i\,,\,v_i+hp_i],
\end{equation}
which are defined by central grid-knots $\msv$. Apparently,
$S_n(\msp,\msv)$ is the union of 
closures of those basic cubes $C_n(\msp,\msx)$
which share the grid-knot $\msv$.

The grids $G_n$ of (\ref{exp3.1}) are homogeneous with respect to
translations in the direction of coordinate axes, i.e 
$\msx \in G_n, \mst = h p_i \mse_i \Rightarrow 
\msx + \mbox{\boldmath{$t $}} \in G_n$
for any $i \in \{1,2,\ldots,d\}$ and $p_i \in {\bbZ}$. There exist
subsets of $G_n$ which are also homogeneous in the defined sense. Let
$\msr_0 \in G_n$ and $\msr =
(r_1, r_2, \ldots, r_d) \in {\bbN}^d$ be fixed. The set
\begin{equation}\label{exp3.2}
 G_n(\msr_0, \msr) \ = \ \{ \msr_0 \:+\: h\,\sum_{l=1}^d \:k_l\,r_l \,
 \mse_l \ : \ k_l \in {\bbZ}\} 
\end{equation}
is a subset of $G_n$ with the following feature $\msx \:\in\: G_n(\msr_0,
\msr), \mst = h p_i r_i\,\mse_i \Rightarrow 
\msx + \mbox{\boldmath{$t $}} \in  G_n(\msr_0, \msr)$.
A grid (\ref{exp3.2}) is denoted by $G_n(R)$, where $R$ stands shortly for 
the $2d$ parameters $\msr_0, \msr$.

Let $h_0 = 2^{-n_0}$ for some $n_0 \in {\bbN}$, $\msp \in {\bbN}^d$ and let
$D$ be a connected set with the structure $D = \cup_{\msvd \in F_n} C_n(\msp,\msv)$,
where $F_n \subset G_n$. For the subgrid $G_n(D) = D \cap G_n(R)$ the
set $G_n(D)$ is discrete and therefore its
interior, closure  and boundary are defined indirectly,
$int\big (G_n(D)\big) = G_n(D) \cap int(D)$, $cls\big (G_n(D)\big)
= G_n(R) \cap \overline{D}$ and $bnd\big (G_n(D)\big)$ is the difference of 
$cls\big (G_n(D)\big)$ and $int\big (G_n(D)\big)$. Apparently,
$int\big (G_n(D)\big) \subseteq G_n(D) \subseteq cls\big (G_n(D)\big)$. 
Let a finite collection of sets $D_l, l \in \EuScript{L}$ make a partition
of ${\bbR}^d$, where each $D_l$ has the structure like the described set $D$.
Then $G(l) = D_l \cap G_n(R)$ make a partition of $G_n$. 

Each $\msx \in G_n$ can be indexed by $\msm
\in {\bbR}^d$, where $\msx = h\msm$.
Similarly, we index grid-knots of
$G_n(\msr_0, \msr)$ by those
$\msm \in {\bbZ}^d$ for which there holds $\msx = \msr_0 + h\sum_l m_l r_l 
\mse_l$. Therefore, we define the sets $I_n = {\bbZ}^d$ and 
$I_n(R) \subset I_n$, indexing the grid-knots of $G_n$ and $G_n(R)$. 
In this work frequently utilized pairs of grids and their index sets are $G_n, 
\ I_n$; $G_n(R), \ I_n(R)$; $G_n(l), \ I_n(l); \ G_n(R,D),I_n(R,D)$.

The {\em shift operator} $Z(\msx), \msx \in {\bbR}^d$, acting on functions 
$f:{\bbR}^d \mapsto {\bbR}$, is defined by $\big(Z(\msx)f\big)(\msx) = 
f(\msx+\msz)$. Similarly we define the discretized shift operator by 
$\big ( Z_n(r,i) {\bf u}_n \big)_{\mskd} = ({\bf u}_n)_{\msld}$, where $\msl = 
\msk + rh\mse_i$.

{\bf Discretization of differential operators}. With respect to a grid step $h$, 
the partial derivatives of $u \in C^{(1)}({\bbR}^d)$ are discretized by 
forward/backward finite difference operators in the usual way,
\begin{equation}\label{ex2.7}\begin{array}{c}
 \der_i(t) u(\msx) \ = \ \frac{1}{t} \big( 
 u(\msx + t \mse_i) \:-\: u(\msx) \big ),\\
 \widehat{\der}_i(t) u(\msx) \ = \ \frac{1}{t} \big( 
 u(\msx) \:-\: u(\msx - t\mse_i) \big ), \end{array} \quad \msx \in {\bbR}^d, \ t \ne 0.
\end{equation}
Let $r \in {\bbZ}\setminus \{0\}$. 
Discretizations of the functions $\partial_iu$ on $G_n$, denoted by
$U_i(r){\bf u}_n, V_i(r){\bf u}_n$, are defined by:
\[ \big (U_i(r)\,{\bf u}_n\big )_{\msmd} \:=\: \der_i\big(rh \big)
 \:u(\msx_{\msmd}), \quad  \big (V_i(r)\,{\bf u}_n\big )_{\msmd} \:=\:
 \widehat{\der}_i \big(rh \big) \:u(\msx_{\msmd}),\]
where $\msx \in G_n$. Then 
\begin{equation}\label{ex2.8}
 \begin{array}{l} U_i(r) \ = \ (rh)^{-1}(Z_n(r,i) \,-\, I),\\
 V_i(r) \ = \  (rh)^{-1} \big(I - Z_n(-r,i) \big) \:= \:
 U_i(-r) \:=\:\,-\, U_i(r)^T.
\end{array}
\end{equation}
Therefore we have $U_i(-r) = U_i(r)\,Z_n(-r,i) = Z_n(-r,i)\,U_i(r)$,
and similarly for $V_i(r)$. In the case of $r=1$ we use a short notations
$U_i, V_i$.

In accordance with the previous terminology, we say that $\partial_i, \sum_{ij}
\partial_i a_{ij} \partial_j$ {\em etc.} are differential operators on ${\bbR}^d$ or $D$.
We say that their discretizations are defined on $G_n$ or $G_n(D)$. In
particular, discretizations of the differential operator (\ref{exp2.1}) are
denoted by $A_n$. Naturally, matrices $A_n$ are the main object in this work.

\subsection{Relations between $l(G_n(R))$ and $W_2^1$-spaces}
The discretization of a function $u \in C({\bbR}^d)$ on $G_n$ is denoted by 
${\bf u}_n$ and defined by values at grid-knots, $\big ({\bf u}_n\big )_{\msmd} \:=\:
u(\msx_{\msmd})$ where $\msx_{\msmd} = (m_1h, m_2h, \ldots,m_dh) \in G_n$, 
and $\msm = (m_1, m_2, \ldots,m_d)$ is a multi-index. The function ${\bf u}_n$ is
usually called a grid-function. We denote the linear spaces of grid-functions 
by $l(G_n)$ or $l(G_n(D))$. The linear space of grid-functions on $G_n(R)$ with 
finite supports is denoted by $l_0(G_n(R))$. Elements of $l(G_n)$ are also
called columns. The corresponding $L_p$-spaces are denoted by $l_p(G_n)$ or 
$l_p(G_n(D))$, and their norms by $\lnorm \cdot \lnorm_p$. The duality pairing
of ${\bf v} \in l_q(G_n)$ and ${\bf u} \in l_p(G_n)$ is denoted by $\lev
{\bf v} |{\bf u} \des$. The scalar product in $l_2(G_n)$ is denoted by
$\lev \cdot |\cdot \des$ and sometimes by $( \cdot |\cdot )$. 
The norm of $l_p(G_n(R))$ is denoted by $\lnorm \cdot \lnorm_{Rp}$. 
For $p \in [1,\infty)$ this norm is defined by:
\[ \lnorm {\bf u} \lnorm_{Rp} \ = \ \Big[ \,\vol(R)\,\sum_{\mskd \in I_n(R)}\:
 |u_{\mskd}|^p\, \Big]^{1/p}, \]
where $\vol(R) = \prod_{i=1}^d r_i$. Finally, for $p = \infty$ we have
$\lnorm {\bf u} \lnorm_{R\infty} = \sup\{|u_{\mskd}|: \msk \in I_n(R)\}$.

Let us define the quadratic functional on $l(G_n)$ by $q({\bf u}) =
\sum_i^d \lnorm U_i{\bf u}\lnorm_2^2$ and
$q_R({\bf u}) = \vol(R)\sum_i^d \lnorm U_i(r_i){\bf u}\lnorm_{R2}^2$ on $l(G_n(R))$.
It is understood $q_R = q$ for $G_n(R) = G_n$. 
A discrete analog of $W_2^1$-spaces is the linear spaces of those 
${\bf u}_n \in l(G_n(R))$ for which the norm $\lnorm \cdot \lnorm_{R2,1}$:
\begin{equation}\label{exp3.3}
 \lnorm {\bf u} \lnorm_{R2, 1}^2 \ = \ \lnorm {\bf u} \lnorm_{R2}^2
 \ + \ q_R({\bf u}),
\end{equation}
is finite. This space is denoted by $w_2^1(G_n(R))$.
By convention $\lnorm \cdot \lnorm_{2,1} = \lnorm \cdot \lnorm_{R2,1}$
for $r_i = 1$. The subspace of grid-functions ${\bf u} \in w_2^1(G_n(R))$ for 
which ${\bf u}_n = {\bbJ}_{G_n(D)}{\bf u}_n$ is denoted by $w_2^1(G_n(R,D))$. 
Hence, $w_2^1(G_n(R,D))$ for $\msr = {\bf 1}$ is denoted by $w_2^1(G_n(D))$. 
For problems on bounded domains we need a discrete version of
the Poincar\'{e} inequality which is formulated as follows:

\begin{lemma}\label{lem3.1} Let $D$ be bounded. Then the norms 
$\lnorm \cdot \lnorm_{2,1}$ and $q_R(\cdot)^{1/2}$ are equivalent in 
$w_2^1(G_n(R,D))$, 
\[ q_R(\cdot)^{1/2} \ \geq \ \beta\,\lnorm \cdot \lnorm_{R2,1},\]
where $\beta$ is independent of $n$.
\end{lemma}

An element (column) ${\bf u}_n \in l(G_n)$ can be associated to a continuous function
on ${\bbR}^d$ in various ways. Here is utilized a mapping $l(G_n) \mapsto C({\bbR}^d)$
which is defined in terms of hat functions. Let $\chi$ be the canonical hat function
on ${\bbR}$, centered at the origin and having the support $[-1,1]$. Then
$z \mapsto \phi(h,x,z) = \chi(h^{-1}(z-hx))$ is
the hat function on ${\bbR}$, centered
at $x \in {\bbR}$ with support $[x-h,x+h]$. The functions
$\msz  \: \mapsto \: \phi_{\mskd}(\msz) \:=\:\prod_{i=1}^d \:\phi(h,x_i,z_i),
x_i = hk_i$, define $d$-dimensional hat functions with supports 
$S_n({\bf 1},\msx) = \prod_i [x_i-h, x_i+h]$. The functions $\phi_{\mskd}(\cdot)
\in G_n$, span a linear space, denoted by $E_n({\bbR}^d)$.
Let ${\bf u}_n \in l(G_n)$ have the entries $u_{n \mskd} = ({\bf u}_n)_{\mskd}$.
Then the function $u(n) =\sum_{\mskd \in I_n} u_{n  \mskd} \phi_{\mskd}$
belongs to $E_n({\bbR}^d)$ and defines imbedding of grid-functions into the
space of continuous functions. We denote the corresponding mapping by
$\Phi_n : l(G_n) \mapsto E_n({\bbR}^d)$. Obviously that there exists
$\Phi_n^{-1} : E_n({\bbR}^d)  \mapsto l(G_n)$ and the spaces 
$l(G_n)$ and $E_n({\bbR}^d)$ are isomorphic with respect to the pair
of mappings $\Phi_n, \Phi_n^{-1}$. It is clear that $E_{n}({\bbR}^d)
\subset E_{n+1}({\bbR}^d)$ and the space of functions $\cup_n E_n({\bbR}^d)$
is dense in $L_p({\bbR}^d), p \in [1,\infty)$, as well as in
$\dot{C}({\bbR}^d)$. Let us mention that $\sum_{\mskd}
\phi_{\mskd} = 1$ on ${\bbR}^d$.

Now we consider another collection of basis functions. To each $\msx = h\msk \in
G_n(R)$ there is associated a $d$-dimensional hat function
\[ \psi_{\mskd}(\msx) \ = \ \prod_{i=1}^d\: \chi \left(\frac{x_i-hk_i}{hr_i}\right ),\]
obviously, with the support $S_n(\msr,\msx) = \prod_i [x_i-r_ih, x_i+r_ih]$. They
span a linear space denoted by $E_n(R,{\bbR}^d)$. Again
we have $\sum_{\mskd}\psi_{\mskd} = 1$ on ${\bbR}^d$.
The mappings $\Phi_n, \Phi_n^{-1}$ cannot be applied to elements
of $l(G_n(R))$ and $E_n(R,{\bbR}^d)$, respectively. Therefore we define
restrictions $\Phi_n(R): l(G_n(R)) \rightarrow E_n(R,{\bbR}^d)$ and the corresponding
inverse mapping $\Phi_n^{-1}(R)$ by the following expression:
\begin{equation}\label{exp3.10}
 u(n) \ = \ \Phi_n(R)\,{\bf u}_n \ = \ \sum_{\mskd}\:
 \big({\bf u}_n\big)_{\mskd}\:\psi_{\mskd}.
\end{equation}
If we have to underline that $u(n)$ is related to a particular set of parameters
$R$ then we use an extended denotation $u(R,n)$. For two functions $v(n), u(n)$ we
have $(v(n)|u(n)) = h^d \vol(R) \sum_{\mskd \msld} s_{\mskd \msld} v_{\mskd} u_{\mskd}$
where $s_{\mskd \msld} = \Vert \psi_{\mskd}\Vert_1^{-1}(\psi_{\mskd}|\psi_{\msld})$.
Let us notice that $\sum_{\msld} s_{\mskd \msld} = 1$.

We cannot compare directly columns ${\bf u}_n$ with various $n$. An
indirect comparison can be made by using $u(n) = \Phi_n(R){\bf u}_n
\in \cup_n E_n(R,{\bbR}^d)$. To compare $U_i(r_i){\bf u}_n$ and $\partial_iu(n)$ we
need an additional expression. Let ${\bf u}$ and $u(n)$ be related by (\ref{exp3.10})
and $\dot{s}_{\mskd\msld}(i) = h^{-d}(\partial_i\psi_{\msk} |\partial_i \psi_{\msl})$. 
Then, for a homogeneous grid $G_n(R), \msr \in {\bbN}^d$, there must hold
\begin{equation}\label{exp3.11} 
 \sum_{\mskd\msld}\: v_{\mskd}  u_{\msld}\, \dot{s}_{\mskd\msld}(i)
 \ = \ -\,\frac{1}{2}\:  \sum_{\mskd\,\msrd'\,r_i}\: 
 \big(v_{\mskd+r_i\msed_i}\,-\, v_{\mskd} \big)\,
 \big(u_{\mskd+\msrd'+r_i\msed_i} \,-\, u_{\mskd+\msrd'}\big)\,
 s_{{\bf 0}\msrd'}'\,\dot{s}_{0r_i},
\end{equation}
where $\msr' = (r_1,r_2,\ldots,r_{i-1},0,r_{i+1},\ldots,r_d)$,
$\dot{s}_{0r}=(\partial_i \psi_0|\partial_i\psi_{r})$ and $s_{{\bf 0}
\msrd'}'$ is the $(d-1)$-dimensional quantity $s_{\mskd\msld}$. 
Then (\ref{exp3.11}) follows from $\sum_{\msld} s_{\mskd \msld} = 1$
and consequently $\sum_{\msld} \dot{s}_{\mskd \msld} = 0$,
after the sum is carried out over any partial component $l_i$ of the index
$\msl$. Thus we have
\begin{equation}\label{exp3.4}
 (\partial_i v(n)\,|\,\partial_i u(n)) \ = \ \Vert \psi_{\mskd}\Vert_1\:
 \sum_{m_i \mskd' \msld'}\: s_{\mskd' \msld'}\, \big(U_i(r_i){\bf v}\big)_{m_i \mskd'}\,
 \big(U_i(r_i){\bf u}\big)_{m_i \msld'},
\end{equation}
where the indices are defined by $\msk' =(k_1,\ldots,k_{i-1},k_{i+1},\ldots,
k_d)$ and analogously $\msl'$.

\begin{lemma}\label{lem3.3} Let sequences of functions $v(n), u(n), 
n \in {\bbN}$, be defined by (\ref{exp3.10}). Then 
\[\begin{array}{lc} (i)&
 \Big| \sum_{i=1}^d \: \big( \,\partial_i v(n)\,|\,\partial_i u(n)\,\big)
 \Big | \ \leq \  h^d\,q_R({\bf v})^{1/2}\,q_R({\bf u})^{1/2}.\\
 (ii)& 
 \Big|\big( \,\partial_i v(n)\,|\,\partial_j u(n)\,\big) \,-\,h^d\,\vol(R)\,
 \sum_{\mskd \in I_n(R)} \:\big(U_i(r_i){\bf v}\big)_{\mskd}\,\big(U_j(r_j)
 {\bf u}\big)_{\mskd} \Big|\\ 
 &\leq \ h^d\,\min\:\left \{ \begin{array}{ll} \displaystyle
 \lnorm U_i(r_i){\bf v}\lnorm_{R2}\,&\sup \:\big\{\lnorm (Z(w,j)-I)\,
 U_j(r_j){\bf u}\lnorm_{R2} \::\: |w| \leq r_jh \big\} 
 \\ \displaystyle
 \lnorm U_j(r_j){\bf u}\lnorm_{R2}\,&\sup \:\big\{\lnorm (Z(w,i)-I)\,
 U_i(r_i){\bf v}\lnorm_{R2} \::\: |w| \leq r_ih \big\}.
 \end{array} \right.  \end{array}\]
\end{lemma}

{\Proof} After applying the CSB-inequality to (\ref{exp3.4}) and using
$\sum_{\msld}s_{\mskd\msld} =1$ we get (i). Assertion (ii) is proved here
for $i=j=1$. In this proof $\partial = \partial_1$. By~(\ref{exp3.4}) 
we can straightforwardly calculate 
\[\begin{array}{c}\displaystyle
 (\partial v(n)\,|\, \partial u(n)) \ = \ \Vert \psi_{\mskd}\Vert_1\:\sum_{\mskd, 
 \msrd'}\: s_{{\bf 0}\msrd'}\, \big(U_i(r_i){\bf v}\big)_{\mskd}\,
 \big(U_i(r_i){\bf u}\big)_{\mskd+\msrd'}\\ \displaystyle
 = \  h^d \,\big(U_i(r_i){\bf v}\big | U_i(r_i){\bf u}\big)_R \:+\:\delta(n),\end{array}\]
where
\[ \delta(n) \ = \ \Vert \psi_{\mskd}\Vert_1\:\sum_{\mskd, 
 \msrd'}\: s_{{\bf 0}\msrd'}\, \big(U_i(r_i){\bf v}\big)_{\mskd}\,\Big[
 \big(U_i(r_i){\bf u}\big)_{\mskd+\msrd'}\,-\, \big(U_i(r_i){\bf u}\big)_{\mskd}
 \Big].\]
By the CBS inequality the error term $\delta(n)$ can be estimated as in Assertion (ii).{\QED}

\begin{lemma}\label{lemn3.4}
Let $G_n(R)$ be a homogeneous subgrid defined by (\ref{exp3.2}).
There exists $\sigma^2 \in (0,1)$ such that
\[ (1\,-\,\sigma^2)\:h^d\:\lnorm {\bf u}_n \lnorm_{R2}^2 \ \leq \
 \Vert \,u(n)\,\Vert_2^2 \ \leq \ h^d\: \lnorm {\bf u}_n \lnorm_{R2}^2\]
uniformly with respect to $n \in {\bbN}$.
\end{lemma}

{\Proof} Let us consider first the one-dimensional case. The grid $G_n(R)$ 
consists of points $x_k = hrk \in {\bbR}, k \in {\bbZ}$, and $E_n(R,{\bbR})$
is spanned by the hat functions $\psi_k$ centered at $x_k$ with the supports
$[-hr+x_k, x_k+hr]$. We define the matrix $S(1)$ with entries:
\[ s_{kl} \ = \ \frac{1}{hr}\,(\psi_k|\psi_l) \ = \ 
 \left \{\begin{array}{lll} (2/3) & {\rm for}& k \:=\:l,\\
 (1/6) & {\rm for}& k \:=\:l\pm r. \end{array} \right .\]
Obviously we have $S(1)=I-(1/3)A$, where the matrix $A$ has the 
structure $A = I+(1/2)(I_++I_-)$ and $I_\pm$ are the first upper and lower
off-diagonals. It is well known \cite{Str} that $A$ has a purely continuous
spectrum in $[0,2]$ so that $S(1)$ has the spectrum equal $[1/3,1]$. Therefore
\[ \Vert u(n) \Vert_2^2 \ = \ h\, \\vol(R)\,\sum_{kl}\:s_{kl}\,u_k u_l
 \geq \frac{1}{3}\, h\, \\vol(R)\:\sum_k \:u_k^2 \ = \ 
 \frac{1}{3}\,h\, \lnorm {\bf u}\lnorm_{R2}^2.\]
Hence, we have here $1-\sigma^2 = 1/3$.

In order to generalize this proof to $d$-dimensional case we proceed as
follows. The symmetric matrix $S(d)$ with entries $s_{\mskd \msld}$
can be represented as the outer product of $d$ matrices $S(1)$
with entries as in the first part of proof. Therefore its spectrum is $Sp(S(d)) =
\prod_{i=1}^d Sp(S(1))$. According to the first part of proof the matrix $S(1)$ has
its spectrum in the interval $[1/3,1]$, implying $Sp(d) \geq 3^{-d}$. Hence,
with $\sigma^2 = 1-3^{-d}$ we have
\begin{equation}\label{exp3.5}
 \vol(R)\:\sum_{\mskd \msld \in I_n(R)}\:s_{\mskd \msld}\,u_{\mskd}\,u_{\msld}
 \ \geq \ (1-\sigma^2)\, \lnorm {\bf u} \lnorm_{R2}^2,
\end{equation}
providing us with a proof of the general case. {\QED}

\begin{theorem}\label{thn3.1} Let $u(n) = \Phi_n(R) {\bf u}_n$. There exists $\sigma^2
\in (0,1)$, independent of $n$, such that
\[\displaystyle
 (1-\sigma^2)\,h^d\,\lnorm {\bf u}_n \lnorm_{R2,1}^2  \ \leq \
 \Vert u(R,n) \Vert_{2,1}^2 \ \leq \ h^d \, \lnorm {\bf u}_n \lnorm_{R2,1}^2. \]
\end{theorem}

{\Proof} It is sufficient to prove the first double inequality. The estimates
from above are obvious. To get the estimates from bellow it suffices to
consider $\partial_i u$. From Expression (\ref{exp3.4}) we have
\[ \Vert \partial_i u(n) \Vert_2^2 \ = \ \Vert \psi_{\mskd}\Vert_1 \:
 \sum_m\: \sum_{\mskd' \msld'} \: s_{\mskd' \msld'}\:
 \big(U_i(r_i){\bf u}\big)_{m \mskd'}\, \big(U_i(r_i){\bf u}\big)_{m \msld'}.\]
Then after applying (\ref{exp3.5}) to the inner sum we get
\[ \Vert \partial u(n) \Vert_2^2 \ \geq \ (1-\sigma^2)\,\Vert \psi_{\mskd}\Vert_1\:
 \sum_{m,\mskd'} \:  \big( U_i(r_i) {\bf u} \big)_{m \mskd'}^2 \ = \
 (1-\sigma^2)\,h^d\:\lnorm U_i(r_i) {\bf u}\lnorm_{R2}^2,\]
from where follows the estimate from bellow. {\QED}

An element $u \in W_2^1({\bbR}^d)$ does not necessary belong to $E_n(R,{\bbR}^d)$.
In order to approximate $u$ with elements of $E_n(R,{\bbR}^d)$ we define:
\begin{equation}\label{exn3.13}
 \hat{u}(n) \ = \ \sum_{\mskd \in I_n(R)}\: \Vert \psi_{\mskd} \Vert_1^{-1}
 \:(\psi_{\mskd}|u)\,\psi_{\mskd}.
\end{equation}
The numbers $\Vert \psi_{\mskd} \Vert_1^{-1}(\psi_{\mskd}|u)$ are called the
{\em Fourier coefficients} of $u$. 

The basic result for our proof of convergence of approximate solutions is
formulated in terms of the quantity
$\Gamma_p(\msw,u)$ defined by:
\[ \Gamma_p(\msw,u) \ = \  \Vert \, (Z(\msw) \,-\,I) \,u\, \Vert_p\]
The kernels
\begin{equation}\label{exp3.6}
 \omega_n(\msx,\msy) \ = \ \sum_{\mskd} \: \frac{1}{\Vert \psi_{\mskd}\Vert_1}\:
 \psi_{\mskd}(\msx)\,\psi_{\mskd}(\msy)
\end{equation}
define an integral operator which is denoted by $K_n$. Actually, the kernels
$\omega_n$ define a $\delta$-sequence of functions on ${\bbR}^d \times {\bbR}^d$
and $K_n$ converge strongly in $L_p$-spaces to unity:
\begin{corollary}\label{corn3.0} Let $p \in [1,\infty]$. Then
\begin{description}\itemsep 0.cm
 \item{(i)} $\Vert K_n \Vert_p \:\leq \:1$.
 \item{(ii)} There is a positive number $\kappa(R)$, independent of $n$, such
that $\Vert (I-K_n)u \Vert_p \:\leq \: \kappa(R)
 \sup\{\Gamma_p(\msw,u)\,:\,|w_i| \leq \:hr_i\}$.
 \item{(iii)} The operator $K_n \in \mathsf{L}(L_2({\bbR}^d),L_2({\bbR}^d))$ has the
spectrum equal $Sp(K_n)=\{0\}\cup [3^{-d},1]$.
\end{description}
\end{corollary}

{\Proof} Only (iii) has to be proved. The symmetric operator $K_n$ is reduced by
$E_n(R,{\bbR}^d)$ and represented by an integral operator with the kernel 
(\ref{exp3.6}).
It is zero operator in the orthogonal complement $E_n(R,{\bbR}^d)^\perp$. With
respect to the mapping $\Phi_n(R)$ the operator $K_n$ is mapped to the symmetric
matrix $\hat{K_n}=\Phi_n(R)^{-1} K_n \Phi_n(R) = S(d)$ in $l_2(I_n(R))$. {\QED}

\begin{theorem}\label{Thn3.2} Let $v,u \in W_2^1({\bbR}^d)$ and 
$\hat{u}(n),\hat{v}(n)$ be defined by (\ref{exn3.13}). Then
\[\begin{array}{l}
 \big |(\hat{v}(n) |\hat{u}(n)) \:-\: (v|u)\big | \ \leq \ 
 c(R)\,\min \:\left \{\begin{array}{l}
 \Vert \,u\,\Vert_2 \:  \sup_{|\msw| \leq h|\msr|}\:\Gamma_2(\msw,v), \\
 \Vert \,v\,\Vert_2 \:  \sup_{|\msw| \leq h\msr|}\:\Gamma_2(\msw,u), \\
 \end{array} 
  \right.
 \\ \\   \displaystyle 
    \big |(\partial_i \hat{v}(n)\,|\,\partial_j \hat{u}(n)) \:-\:
 (\,\partial_i v\,|\, \partial_j u\,)\big | \\
  \leq \  c(R)\,\min\:\left \{ \begin{array}{l} \displaystyle
 \Vert \partial_i v \Vert_2\,\Big [  \Vert \partial_j u - 
 \der_j u\Vert_2 \,+\, \sup_{|\msw| \leq h\msr|}  \: \Gamma_2(\msw,\partial_j u) \Big ],
 \\ \displaystyle
 \Vert \partial_j u \Vert_2\,\Big [  \Vert \partial_i v - 
 \der_i v\Vert_2 \,+\, \sup_{|\msw| \leq h\msr|}  \: \Gamma_2 (\msw,\partial_i v) \Big ], 
 \end{array} \right. \end{array}
  \]
where $c(R)$ is $n$-independent.
\end{theorem}

A proof of this theorem is rather technical. For instance, in order to prove the
second inequality one has to use a sequence of replacements: 
$(\partial_i \hat{v}(n)|\partial_j \hat{u}(n)) \rightarrow
(\der_i(h)v|K_n\der_j(h)u) \rightarrow (\der_i(h)v|K_n\partial_ju) \rightarrow
(\partial_i v|K_n\partial_j u) \rightarrow (\partial_i v|\partial_j u)$. Each
replacement gives rise to an error. The sum of errors can be estimated by an
expression as given in the second inequality of assertion.

\subsection{Imbedding of $l(G_n(R))$ into $W_2^{-1}$-spaces}\label{sec32}
Beside the functions $\psi_{\mskd}$ we consider the functions defined by:
\[ \chi_{\mskd i +}(\msx) \ = \ {\bbJ}_{[k_i,k_i+ r_ih)}(x_i)\:
 \prod_{j \neq i}\: \psi_{k_j}(x_j),\quad \chi_{\mskd i -}(\msx) \ = \ 
 \chi_{\mskd -r_i\msed_i i +}(\msx) \]
for all the possible $i = 1,2,\ldots,d$, and the linear space $F_n(R,{\bbR}^d)$
spanned by the defined functions $\chi_{\mskd i \pm}$. 
The integral operators $K_n^{(i)}(\chi,\chi)$ with the respective kernels
\[ \omega_n^{(i)}(\msx,\msy) \ = \ \sum_{\mskd} \: \frac{1}{\Vert \chi_{\mskd i +}\Vert_1}\:
 \chi_{\mskd i +}(\msx)\,\chi_{\mskd i +}(\msy)\:=\:
 \sum_{\mskd} \: \frac{1}{\Vert \chi_{\mskd i +}\Vert_1}\:
 \chi_{\mskd i -}(\msx)\,\chi_{\mskd i -}(\msy)\]
have properties similar to the integral operators $K_n$ of Corollary \ref{corn3.0}.
The same is valid for the non-symmetric integral operators $K_n^{(i)}(\psi,\chi)$,
and their adjoints $K_n^{(i)}(\chi,\psi)$, where the kernel of $K_n^{(i)}(\psi,\chi)$ is
defined by $(\msx, \msy) \mapsto (h^d \vol(R))^{-1}\sum_{\mskd}\psi_{\mskd}(\msx)
\msch_{\mskd i +}(\msy)$.

\begin{lemma}\label{lem3.6} The operators $K_n^{(i)}(\chi,\chi), K_n^{(i)}(\psi,\chi),
K_n^{(i)}(\psi,\chi)^\dag$ have properties (i) and 
(ii) of Corollary \ref{corn3.0}. The spectra of operators $K_n^{(i)}(\chi,\chi) \in 
\mathsf{L}(L_2({\bbR}^d),L_2({\bbR}^d))$ consist of $0$ and an interval $[\kappa,1]$
with certain $\kappa \in (0,1)$.
\end{lemma}

An element $\mu \in W_2^{-1}({\bbR}^d)$ is represented as $\mu = f_0+\sum_{i=1}^d
\partial_i f_i$ with $f_i \in L_2({\bbR}^d)$ with the norm $\Vert \mu \Vert_{2,-1}
= (\sum_i \Vert f_i\Vert_2^2)^{1/2}$. Its discretizations are defined by
grid-functions $\msmu_n$ with the components:
\begin{equation}\label{exp3.7}
 \mu_{\mskd} \ = \ \frac{1}{\Vert \psi_{\mskd}\Vert_1}\,(\psi_{\mskd}|f_0)
 \,-\,\frac{1}{\Vert \psi_{\mskd}\Vert_1}\:\sum_{i=1}^d\:(\partial_i
 \,\psi_{\mskd}| f_i).
\end{equation}
Therefore we can write:
\begin{equation}\label{exp3.8}
 \lev u(n) \,|\,\mu\des \ = \ h^d\:\lev {\bf u}_n\,|\,\msmu_n\des_R \ = \
 \Vert \psi_{\mskd}\Vert_1\:\sum_{\mskd \in I_n(R)}\:u_{\mskd}\mu_{\mskd}.
\end{equation}
The expression
\[ \partial_i \psi_{\mskd} \ = \ \frac{1}{r_ih}\:\Big[\chi_{\mskd i -}-
 \chi_{\mskd i +}\Big], \quad a.e. {\rm~~on~~} {\bbR}^d,\]
enables us to rewrite the components of $\msmu_n$ in (\ref{exp3.7}) in terms of the
Fourier coefficients $\hat{f}_{i \mskd \pm}=
\Vert \psi_{\mskd}\Vert_1^{-1}(\chi_{\mskd i \pm}|f_i)$ of functions $f_i$.
The corresponding grid-functions $\hat{{\bf f}}_{in \pm}$ are imbedded into
the space $F_n(R,{\bbR}^d)$ by the mappings $\check{{\bf f}}_{in \pm} \mapsto 
\check{f}_{i \pm}(n)= K_n^{(i)}(\chi,\chi)f_i$ which are completely analogous to the
mapping $\hat{{\bf f}}_0 \mapsto \hat{f}_0(n) = K_n f_0$ of Corollary
\ref{corn3.0}. Now we can get the following useful
expression:
\begin{equation}\label{exp3.9}
 \lev u(n) \,|\,\mu\des \ = \ h^d\:\Big[ ({\bf u}_n|\hat{{\bf f}}_{0n})_R
 \,-\,\sum_{i=1}^d\: (U_i(r_i){\bf u}_n|\check{{\bf f}}_{in +})\Big].
\end{equation}
Therefore, by~(\ref{ex2.8}) we can represent elements $\msmu \in l(G_n(R))$ 
in the following way:
\begin{equation}\label{exp3.12}
 \msmu_n \ = \ \hat{{\bf f}}_{0n}\:-\:\sum_{i=1}^d\:V_i(r_i)\check{{\bf f}}_{in+}.
\end{equation}
From Expressions (\ref{exp3.8}) and (\ref{exp3.9}) we have $|\lev {\bf u}_n|
\msmu_n\des| \leq \lnorm {\bf u}_n\lnorm_{2,1} (\lnorm \hat{{\bf f}}_0\lnorm_2^2
+\lnorm \check{{\bf f}}_{in+}\lnorm_2^2)^{1/2}$, implying a natural definition:
\[ \lnorm \msmu_n \lnorm_{2,-1}^2 \ = \ \sum_{i=1}^d\,\lnorm \hat{{\bf f}}_0
 \lnorm_2^2 +\lnorm \check{{\bf f}}_{in+}\lnorm_2^2.\]

For the sake of a concise writing of final results we denote here $K_n$ by 
$K_n^{(0)}$.

\begin{lemma}\label{lem3.5} For each $\mu \in W_2^{-1}({\bbR}^d)$ and the
corresponding discretizations $\msmu_n$ defined by (\ref{exp3.12}) the following
is valid:
\[ h^d\,\lnorm \msmu_n \lnorm_{2,-1}^2 \ = \ \sum_{i=0}^d \:(f_i\,|\,K_n^{(i)}
 (\chi,\chi)f_i) \ \leq \ \Vert \mu \Vert_{2,-1}^2,\]
and
\[ 0 \ \leq \ \Vert \mu \Vert_{2,-1}^2\:-\:h^d\,\lnorm \msmu_n \lnorm_{2,-1}^2 \ 
 = \ \sum_{i=0}^d \:(f_i\,|\,(I-K_n^{(i)}(\chi,\chi))f_i) \ \to \ 0,\]
as $n \to \infty$.
\end{lemma}

\subsection{Imbedding into $W_p^1$-spaces}
Theorem \ref{thn3.1} can be partially generalized.

\begin{lemma}\label{lem3.7} The following assertions are valid:
\begin{description}\itemsep 0.cm
 \item{(i)} If $p \in [1,\infty]$ and ${\bf u}_n \in l_p(G_n(R))$ then:
\[ \Vert \,u(n)\,\Vert_{p,1} \ \leq \ h^{d/p}\,\lnorm {\bf u}_n \lnorm_{Rp,1}.\]
 \item{(ii)} If $p \in [1,\infty]$ and $u \in W_p^{1}({\bbR}^d)$ then:
\[ \Vert \hat{u}(n)\Vert_{p,1} \ \leq \ 
 h^{d/p}\,\lnorm \hat{{\bf u}}_n \lnorm_{Rp,1} \ \leq \ 
 \Vert u \Vert_{p,1}.\]
 \item{(iii)} Let $u = u^++u^-$, where $u^+ = \max\{u,0\}, u^-=\min\{u,0\}$. Then
\[ h^d \,\lnorm \hat{{\bf u}}_n \lnorm_{R1} \ \leq \ \Vert \,u\,\Vert_1\:
 =\:\Vert \hat{u}^+(n)\Vert_{1}+\Vert \hat{u}^-(n)\Vert_1 \ \leq \ h^d\,
 \big(\lnorm \hat{{\bf u}}_n^+ \lnorm_{R1}\,+\,\lnorm \hat{{\bf u}}_n^-
 \lnorm_{R1} \big).\]
\end{description}
\end{lemma}

{\Proof} Let us consider the mapping $\Phi_n(R): l(G_n(R)) \mapsto E_n(R,{\bbR}^d)$
defined by $u(n) = \sum_{\mskd} u_{\mskd}\psi_{\mskd}$. In terms of functions 
$\chi_{\mskd i +}$ of Subsection \ref{sec32} we easily get:
\begin{equation}\label{exp6.8}
 \partial_i \,u(n) \ = \ \sum_{\mskd \in I_n(R)}\: \big(U_i(p_i){\bf u}_n
 \big)_{\mskd}\,\chi_{\mskd i +}.
\end{equation}
Apparently we have for $p = 1$ and $p=\infty$ the following inequalities: $\Vert u \Vert_p
\leq h^{d/p}\lnorm {\bf u}_n\lnorm_{Rp}$, $\Vert \partial_i u \Vert_p \leq h^{d/p}\lnorm 
U_i(p_i){\bf u}_n\lnorm_{Rp}$. Hence, by the Riesz-Thorrin theorem we get (i).

Similarly we prove (ii). For $p=1$ and $p=\infty$ we have $h^{d/p}\lnorm 
\hat{{\bf u}}_n \lnorm_{Rp}  \leq \Vert u\Vert_p$. Now we use the expression
$\psi_{\mskd+p \msed_i}(\msx) = \psi_{\mskd}(\msx-ph \mse_i)$ and get
\[ \big(U_i(p) \hat{{\bf u}}_n\big)_{\mskd} \ = \ \frac{1}{\Vert \psi_{\mskd}
 \Vert_1}\: (\psi_{\mskd}\,|\, \der_i(p)\,u).\]
Because of $\Vert \der_i(p)\,u\Vert_p \leq \Vert \partial_i u\Vert_p$ for
$p = 1,\infty$, we have $h^{d/p}\lnorm U_i(p) \hat{{\bf u}}_n \lnorm_{Rp} \leq 
\Vert \partial_i u\Vert_p$. Assertion (ii) follows now from the Riesz-Thorrin theorem. 

The right hand side of double inequality (iii) is implied by (i) while the
left hand side is implied by (ii). Therefore we
have to prove the equality $\Vert u\Vert_1 = \Vert u^+(n)\Vert_{1}+\Vert u^-(n)
\Vert_1$. This equality is a consequence of $\sum_{\mskd}\psi_{\mskd}=1$ on ${\bbR}^d$ 
and $\Vert\hat{u^+}(n)\Vert_1 = \sum_{\mskd}(\psi_{\mskd}|u^+)$.
{\QED}

\section{Construction of discretizations}\label{sec4}
It is important to underline at the beginning that discretizations $A_n$ of differential 
operator $A(\msx)$ are defined prior to discretizations of the forms (\ref{exp2.9}),
(\ref{exp2.10}). This fact is in a full agreement with construction of discretizations 
$A_n$ in this section. Some classes of discretizations are derived from a
general principle which is not based on finite difference formulas and cannot be
apriori related to variational equalities. Nevertheless, 
bilinear forms must be associated to $A_n$ so that $A_n$ are derived from the
corresponding variational equalities. The constructed bilinear forms are considered
as discretizations of the original form (\ref{exp2.4}). These forms are basic objects in our 
proof of convergence of approximate solutions. In the next two subsections schemes
and the corresponding discretized forms are constructed for two classes of methods.

Forms $a_n(\cdot,\cdot)$ on $E(R,{\bbR}^d) \times E(R,{\bbR}^d)$ and
matrices $A_n$ on $G_n(R)$ are related by equalities:
\[ a_n(v,u) \ = \ \lev {\bf v}\,|\,A_n\,{\bf u}\des_R. \]
In addition, the discretized forms determine the discretized variational equalities:
\begin{equation}\label{exp4.1}
 \lambda \lev {\bf v}_n\,|\,{\bf u}_n\des_R \:+\:a_n(v,u) \ = \
 \lev{\bf v}_n\,|\,\msmu_n\des_R, \quad {\bf v}_n \,\in\, w_{2,1}(G_n(R)).
\end{equation}

To discretize $A(\msx)$ means to associate to $A(\msx)$ a sequence of matrices
$A_n$ on $G_n(R), n \in {\bbN}$. Of course, the matrices $A_n$ must be constructed
reasonably in order to enable demonstrations of the convergence of numerical
solutions. The convergence analysis is postponed until two next sections. Therefore,
in this section, the terminology "discretizations" of $A(\msx)$ instead of
approximations of $A(\msx)$ seams to be more suitable.

Discretizations to be considered in this section are possible if certain conditions 
on $a_{ij}$ are fulfilled. The required conditions are stronger than in 
Assumption \ref{Ass2.1}. By relaxing them gradually as $n \to \infty$ we obtain
discretizations for a general $A(\msx)$ given by Assumption \ref{Ass2.1}.

To a given diffusion tensor $a = \{a_{ij}\}_{11}^{dd}$ we associate an auxiliary
tensor $\hat{a}$ defined by the expressions:
\begin{equation}\label{expaa}
 \hat{a}_{ii} \ = \ a_{ii}, \quad \hat{a}_{ij} \ = \ -|a_{ij}| \quad i \ne j.
\end{equation} 

\begin{assumption}[on discretization conditions]\label{ass4.1}
 \item{1.} There exist $\msq \in {\bbN}^d$, a finite index set $\EuScript{L}$ and a partition 
${\bbR}^d = \cup_l D_l$, where the sets $D_l$ are connected unions of $C_n(\msq,
\msx), \msx \in G_n(Q)$. The sets $D_l$ and the diffusion tensor $a = \{a_{ij}\}_{11}^{dd}$
must be mutually related as follows:
\begin{itemize}\itemsep 0.cm
 \item[-] The tensor-valued function $\msx \mapsto a(\msx)$ is continuous on $D_l$
and the functions $a_{ij}, i \ne j$ do not change signs on $D_l$.
 \item[-] There is a grid-step $h_\varepsilon = 2^{-n(\varepsilon)}$ and the closed sets
\[ D_l(\varepsilon) \ = \ \cup_{\msxd \in cls(D_l)}\: S_{n(\varepsilon)}(\msq,
 \msx),\]
such that the functions $a_{ij}$ can be extended to $D_l(\varepsilon)$, not changing 
the signs on $D_l(\varepsilon)$, and the strict ellipticity (\ref{exp2.2}) is valid on 
$D_l(\varepsilon)$ with the same bounds $\Md, \Mg$. 
\end{itemize}
 \item{2.} For each $l \in \EuScript{L}$ the auxiliary diffusion tensor 
$\hat{a}$ is strictly positive definite on $D_l(\varepsilon)$.
 \item{3.} To each $D_l$ there is associated parameter $\msp(l) \in {\bbN}^d$, such that
the following inequality is valid:
\[\begin{array}{ll} \displaystyle
 \omega(a) \ = & \inf_n\: \min_{l \ i} \displaystyle\inf \Big\{\,
 \frac{1}{p_i(l)}\, \inf_{\mszd \in S_n(\mspd(l),\msxd)}\,a_{ii}(\msz)\\ & \displaystyle
 -\: \sum_{m \ne i} \,\frac{1}{p_m(l)}\, \sup_{\mszd \in S_n(\mspd(l),\msxd)}
 \,|a_{im}(\msz)| \::\: \msx \in G_n \cap D_l(\varepsilon)  \Big\} \ > 0.
 \end{array}\]
\end{assumption}

Condition 3. is crucial in our construction of discretizations $A_n$ which
have a particular feature called the compartmental structure. In the next
definition $I$ is an index set, and matrices $A=\{a_{ij}\}_{II}$ are
considered in linear spaces $l_p(I)$ consisting of functions on the set $I$:

\begin{definition}[Compartmental structure]\label{def4.1}
A matrix $A$ in $l_\infty(I)$ is said to be of positive type if 
$A = pI-B$, $p > 0, \ B \geq 0$ and $\lnorm B\lnorm_\infty \leq p$. 
It is called conservative if $B {\bf 1} = p{\bf 1}$.
A matrix $A$ in $l_1(I)$ is said to have the compartmental
structure if $A = pI-B$, $B \geq 0$ and $\lnorm B\lnorm_1 \leq p$. It is
called conservative if for each ${\bf u} \geq {\bf 0}$ there holds $\lnorm B 
{\bf u}\lnorm_1 = p\lnorm {\bf u} \lnorm_1$. 
\end{definition}

\begin{lemma}\label{lem4.1} Let $A$ in $l_1(I)$ be compartmental. Then
 \item{(i)} If $A$ is conservative then $\lnorm B^m \lnorm_1 = p^m$ for each 
$m \in {\bbN}$.
 \item {(ii)} The spectrum $sp(A)$ of a compartmental matrix is contained in $\Re 
\lambda \geq 0$. If $A$ is conservative, then $0 \in sp(A)$.
\end{lemma}

In the case of $d = 2$ Conditions 2. and 3. can be always fulfilled.
For $d > 2$ there exist positive definite diffusion tensors 
$\{a_{ij}\}_{11}^{dd}$ such that Conditions 2. and 3. are not possible \cite{LR3}. 

Discretizations $A_n$ are defined in terms of its matrix entries $(A_n)_{\mskd\msld}$,
where $h\msk, h\msl \in G_n$. For a fixed $\msx = h\msk \in G_n$ the set of all
the grid-knots $\msy = h\msl$ such that $(A_n)_{\mskd\msld} \ne 0$ is denoted by
${\cal N}(\msx)$ and called the {\em numerical neighbourhood} of $A_n$ at $\msx \in G_n$.
In our constructions the sets ${\cal N}(\msx)$ for $\msx \in G_n \cap D_l$
are mutually alike. A set ${\cal N}(\msx)$ contains always a "cross" consisting
of $\msx$ and $2d$ elements $\pm hp_i(l)\mse_i$. Additional elements of
${\cal N}(\msx)$ depend on the sign of $a_{ij}, i \ne j$.

Two classes of discretizations are analyzed.
One of these classes can be defined straightforwardly in terms
of forward and backward difference formulas. The resulting discretizations 
are called {\em basic schemes}. Discretizations of the other class come from a 
general principle \cite{LR3} and they are called {\em extended schemes}.

The compartmental structure of discretizations $A_n$ of the differential operator
$A(\msx)$ is the goal of overall analysis. Here
we describe a general approach to the constructions of discretizations $A_n$
with the compartmental structure which is based on reduction to the 
two-dimensional problems.

The index set of pairs $I(d) = \{ \{ij\} : i < j,\:
i,j = 1,2,\ldots,d, i \ne j\}$ has the cardinal number $m(d)=d(d-1)/2$.
To each index $\{kl\} \in I(d)$ we associate three coefficients,
\begin{equation}\label{exp4.12}
 a_{kk}^{\{kl\}} = \frac{1}{d-1}\, a_{kk}, \quad a_{ll}^{\{kl\}} =
 \frac{1}{d-1}\, a_{ll}, \quad a_{kl}^{\{kl\}} = a_{kl},
\end{equation}
and a bilinear form $a^{\{kl\}}(\cdot, \cdot)$,
\begin{equation}\label{exp4.13}
 a^{\{kl\}}(v,u) \ = \sum_{i,j \in \{r,s\}}\: \int_D \: a_{ij}^{\{kl\}}(\msx) \:
 \partial_i v(\msx)\, \partial_j u(\msx) \,d\msx.
\end{equation}
Apparently, for each pair $v, u \in C^{(1)}({\bbR}^d)$ 
with compact supports, the following equality is valid:
\[ a(v,u) \ = \ \sum_{\{kl\} \in I(d)} \, a^{\{kl\}}(v,u) . \]
To each of the forms $a^{\{kl\}}(\cdot, \cdot)$ we must associate a sequence of
forms $a_n^{\{kl\}}(\cdot,\cdot)$ and matrices $A_n^{\{kl\}}$ constructed by 
two-dimensional schemes. Then the matrix 
\begin{equation}\label{exp4.14}
 A_n \ = \ \sum_{\{kl\} \in I} \:A_n^{\{kl\}} ,
\end{equation}
is a discretization of $A_0(\msx)$. If each $A_n^{\{kl\}}$ has the 
compartmental structure then $A_n$ also has the compartmental structure.
However, $A_n$ can have the compartmental structure although no $A_n^{\{kl\}}$ 
is compartmental. Condition 3. of Assumption \ref{ass4.1} ensures this
advantageous property in our constructions.

\subsection{Two methods of discretizations} 
The forms must be constructed 
by the rules (\ref{exp4.12})-(\ref{exp4.13}) having in mind that
the construction for higher dimensional cases depends on the construction for
two-dimensional case. Therefore, we are due to describe the construction for the
two-dimensional case. 

\subsection*{Basic schemes}
In the case of numerical grids $G_n(P)$ the numerical schemes and corresponding
discrete bilinear forms can be easily mutually related. A proof of convergence
in Sobolev spaces for various right hand sides are presented in our works 
\cite{LR2,LR3}. In the present case we extend analysis to problems with the 
function $a_{12}$ having both signs. Difficulties appear at those grid-knots 
where the sign of $a_{12}$ changes.
\begin{figure}
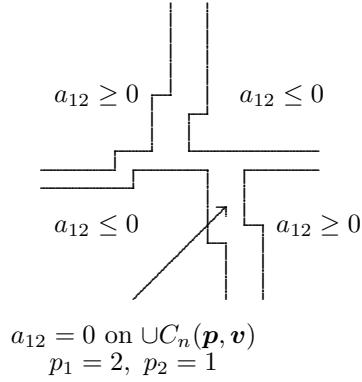

\caption{Asumption on function $a_{12}$}
\label{fig4.0}
\beginpicture
\setcoordinatesystem units <0.7pt,0.7pt>
\setplotarea x from -200 to 0, y from 10 to 220
\setlinear
\plot 20 100, 70 100, 70 110, 110 110, 110 70, 120 70, 120 40 /
\plot 140 40, 140 80, 130 80, 130 110, 170 110 /
\plot 20 110, 60 110, 60 120, 80 120, 80 150, 90 150, 90 200 /
\plot 110 200, 110 140, 100 140, 100 120, 170 120 /
\plot 70 40, 120 90 /
\plot 115 90, 120 90, 120 85 /
\put{$p_1 = 2, \ p_2 =1$} [c] at 70 5
\put{$a_{12}=0$ on $\cup C_n(\msp,\msv)$} [c] at 70 20
\put{$a_{12} \leq 0$} [c] at 50 80
\put{$a_{12} \leq 0$} [c] at 150 150
\put{$a_{12} \geq 0$} [c] at 50 150 
\put{$a_{12} \geq 0$} [c] at 170 80
\endpicture
\end{figure}
\par
The discretizations are constructed by assuming that the sets 
$\{\msx \in {\bbR}^d:a_{12}(\msx) < 0\}$ and $\{\msx \in {\bbR}^d:a_{12}(\msx) > 0\}$
are separated by a set (connected) so that it contains a connected subset
which is equal to the union of cubes $C_n(\msp,\msx), \msx \in G_n(P)$ for some
$n \in {\bbN}$. An illustration of this assumption is given in Figure \ref{fig4.0}. 
This assumption is not valid for a general matrix-valued function $\msx \mapsto
a(\msx)$. Therefore we have to comprehend this assumption as a step of
an approximation procedure in our process of construction of discretizations
$A_n$. For the sake of simple and brief presentation we assume in the next
construction that the assumption is valid already for $n=1$.

Let us define the sets $D_n(-) \subset {\bbR}^d$ by:
\[ D_n(-) \ = \ \{\cup_{\msvd \in G_n(P)} C_n(\msp,\msv) \,:\, a_{12}(\msx) 
 \leq 0 \ {\rm for} \ \msx \in C_n(\msp,\msv)\},\]
and $D_n(+) = D_n(-)^c$. Now we define the subgrids $G_n(P,-)$ consisting of 
all the vertices $\msv \in G_n(P)$ which determine the set $D_n(-)$. As well
we need the subgrid $G_n(P,+)$ consisting of those grid-knots $\msx$ for
which the segments $I(\msx,\msy), \msy = \msx +h_np_1\mse_1$ have the
following property $I(\msx,\msy) \subset cls(D_n(+))$. It is easy to
verify that each $\msx \in G_n(P)$ must be contained in one of sets $G_n(P,\mp)$
and each of segments $I_i(\msx,\msy), \msy=\msx+hp_i\mse_i$ must have both of
its end points $\msx, \msy$ in some $G_n(P,\mp)$.
Some of grid-knots and some of segments are contained in both sets,
$G_n(P,\mp)$ and $cls(D_n(\mp))$, respectively. 
For a segment in $\mse_2$-direction, $I(\msx_1,\msx_2), \msx_2=\msx_1+h_np_2
\mse_2$, the following is true. If $\msx_1 \in G_n(P,-)$ then $I(\msx_1,\msx_2)
\subset cls(D_n(-))$. If $\msx_2 \in G_n(P,+)$ then $I(\msx_2,\msx_1)$ (downward
vertical segment) may be outside of $cls(D_n(+)))$. If this happens then
this segment is contained in a cube $cls(C_n(\msp,\msz))$ on which the function
$a_{12}$ has zero values. This fact is a consequence of a strict separation
of supports of functions $\max\{a_{12},0\}$ and $\min\{a_{12},0\}$ and
will be utilized in our constructions of discrete bilinear forms.

To define forms and entries we need
\[\begin{array}{l}
 \msz^{(\pm)} \ = \ p_i\mse_i \,\pm\,p_j\mse_j \ \in I_n,
 \\
 \msx^{(\pm +)}(n) \ = \ \frac{1}{2}\,\big (\pm h\,p_i
 \mse_i \,+\,h\,p_j\mse_j\big)
 \ \in S_n(\msp,{\bf 0}),\\
 \msx^{(\pm -)}(n) \ = \ \frac{1}{2}\,\big (\pm h\,p_i
 \mse_i \,-\,h\,p_j\mse_j\big)
 \ \in S_n(\msp,{\bf 0}).\end{array}\]
Obviously we have $h\msz^{(\alpha \beta)} \in G_n$ while $\msx^{(\pm)}(n)$ 
are not necessary in $G_n$.

Let us consider a sequence of two-dimensional forms on $E_n(R,{\bbR}^2) \times
E_n(R,{\bbR}^2)$ which are defined by the following expressions:
\begin{equation}\label{exp4.5}
\begin{array}{lll}
 a_n(v,u) &=& a_n^{(-)}(v,u)\:+\: a_n^{(+)}(v,u)\\ \displaystyle
 a_n^{(-)}(v,u) &=& \sum_{i\,j \,=1}^2\:
 \sum_{\msxd \in G_n(P,-)}\: \big(\der_i(p_ih)v\big)(\msx)
 \\ &\times&\displaystyle
 a_{ij}(\msx+\msx^{(++)}
 (n)) \: \big(\der_j(p_jh) u\big)  (\msx),\\ \displaystyle
 a_n^{(+)}(v,u) &=& \sum_{i\,j \,=1}^2\:
 \sum_{\msxd \in G_n(P,+)}\:\big(\der_i((-1)^{i-1}p_ih)v\big)(\msx) 
 \\ &\times&\displaystyle
 a_{ij}(\msx+\msx^{(+-)}(n)) \: \big( \der_j((-1)^{j-1}
 p_jh) u\big) (\msx). \end{array}
\end{equation}
Discretizations $A_n$ of differential operator $A_0(\msx)$ can be easily
obtained from the constructed forms variationaly.
Let us define matrices 
$\EuScript{A}_n^{(\pm)}(i,j,p_i,p_j)$ as the diagonal matrices with entries
$\EuScript{A}_n^{(\pm)}(i,j,p_i,p_j)_{\msxd \msxd}=a_{ij}(\msx+
\msx^{(+\pm)}(n))$. In terms of matrices $U_i(p_i), V_i(p_i)$ 
and $Z_n(p_i,i)$ of (\ref{ex2.8}) we get the following expressions:
\begin{equation}\label{exp4.19}\begin{array}{ll}
 A_n \ =& \displaystyle
 -\:\sum_{l \in \EuScript{L}_-}\:\sum_{ij =1}^2\:V_i(p_i)
 \,\EuScript{A}_n^{(+)}(i,j,p_i,p_j)\,U_j(p_j)\,{\bbJ}_{G_n(P,-)} \\
 & \displaystyle
 -\: \sum_{l \in \EuScript{L}_+}\:\sum_{ij =1}^2\:V_i(p_i)\,
 \Lambda_i^T\,\EuScript{A}_n^{(-)}(i,j,p_i,p_j)\,\Lambda_j
 \,U_j(p_j)\,{\bbJ}_{G_n(P,+)},\end{array}
\end{equation}
where $\Lambda_1 = I, \Lambda_2 =Z_n(-p_2,2)$ and ${\bbJ}_{G_n(P,\mp)}$ are
the projectors on the linear subspaces of grid-functions with
supports in the sets $G_n(P,\mp)$, respectively. The entries of $A_n$ can
be also easily calculated for grid-knots of $int(G_n(P,\mp))$. 
In order to get simple expressions we use the following abbreviations:
\[ a_{ij}^{(\alpha \beta)} \ = \ a_{ij}( \msx + \msx^{(\alpha \beta)}(n)),
 \quad \alpha, \beta \in \{+,-\}.\]
The entries on the cross $\msx \pm h_rp_r\mse_r, r = 1,2$ 
have the structure:
\[\begin{array}{ll}\displaystyle
 \big(A_n\big)_{\mskd \mskd\pm p_1\msed_1}& \displaystyle
 = \  -\,\frac{1}{h^2p_1^2} 
 \: \left\{\begin{array}{lll} \displaystyle
 a_{11}^{(\pm +)}\:-\: \frac{p_1}{p_2}\,
 \big|a_{12}^{(\pm +)}\big| &{\rm for}& a_{12} \leq 0, \\ 
 \displaystyle
 a_{11}^{(\pm -)}\:-\: \frac{p_1}{p_2}\,
 \big|a_{12}^{(\pm -)}\big| &{\rm for}& a_{12} \geq 0,  \end{array}
 \right . \\ \displaystyle

 \big(A_n\big)_{\mskd \mskd\pm p_2\msed_2} & \displaystyle
 = \  -\,\frac{1}{h^2p_2^2} 
 \: \Big[ a_{22}^{(+\pm)}\:-\: \frac{p_2}{p_1}\,
 \big|a_{12}^{(+\pm)}\big|\Big].  \end{array}\]
The entries in the plane spanned by $\mse_1, \mse_2$ have the structure:
\[\begin{array}{l}\displaystyle
 \big(A_n\big)_{\mskd \mskd \pm \mszd^{(-)}} \ = \  
 -\,\frac{1}{h^2p_1p_2}\:
 \big|a_{12}^{(\pm \mp)}\big| \quad {\rm for}\quad a_{12} \leq 0 \ {\rm on} \ D_l,
 \\ \displaystyle
 \big(A_n\big)_{\mskd \mskd \pm \mszd^{(+)}} \ = \  
 -\,\frac{1}{h^2p_1p_2}\:
 \big|a_{12}^{(\pm \pm)}\big| \quad {\rm for} \quad a_{12} \geq 0 \ {\rm on} \ D_l,
 \end{array}\]
where $+\msz^{(\pm)}$ is associated with the upper and $-\msz^{(\pm)}$ with the
lower indices of $a^{(\alpha \beta)}$. The diagonal entries $(A_n)_{\mskd \mskd}$
are equal to the negative sum of all the entries $(A_n)_{\mskd \msld}$,
$\msl \ne \msk$. Entries for grid-knots at $bnd(G_n(P,\mp))$ can be more
complex.
\begin{figure}
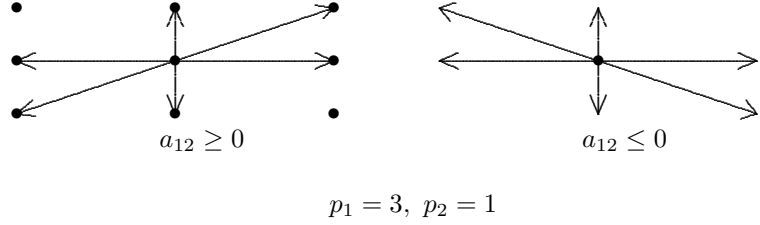

\caption{Numerical neighbourhoods at two grid-knots of $G_n(P,\mp)$}
\label{fig4.1}
\beginpicture
\setcoordinatesystem units <1.0pt,1.0pt>
\setplotarea x from -40 to 300, y from 0 to 100
\setlinear
\plot 0 60, 120 60 /
\plot 60 40, 60 80 /
\plot 0 40, 120 80 /
\plot 160 60, 280 60 /
\plot 220 40, 220 80 /
\plot 160 80, 280 40 /
\plot 7 57, 0 60, 7 63 /
\plot 113 57, 120 60, 113 63 /
\plot 7 39, 0 40, 6 45 /
\plot 116 76, 120 80, 113 81 /
\plot 57 47, 60 40, 63 47 /
\plot 57 73, 60 80, 63 73 /
\plot 167 57, 160 60, 167 63 /
\plot 273 57, 280 60, 273 63 /
\plot 217 47, 220 40, 223 47 /
\plot 217 73, 220 80, 223 73 /
\plot 166 82, 160 80, 164 75 /
\plot 274 38, 280 40, 276 45 /
\put{$\bullet$} [c] at 60 60
\put{$\bullet$} [c] at 0 60
\put{$\bullet$} [c] at 120 60
\put{$\bullet$} [c] at 0 80
\put{$\bullet$} [c] at 60 80
\put{$\bullet$} [c] at 120 80
\put{$\bullet$} [c] at 0 40
\put{$\bullet$} [c] at 60 40
\put{$\bullet$} [c] at 120 40
\put{$\bullet$} [c] at 220 60
\put{$p_1 = 3, \ p_2 =1$} [c] at 150 5
\put{$a_{12} \geq 0$} [c] at 70 30
\put{$a_{12} \leq 0$} [c] at 230 30
\endpicture
\end{figure}
In two-dimensional problems the forms (\ref{exp4.5}) are natural and the 
entries of $A_n$ calculated from variational equalities (\ref{exp4.1}) 
can be perceived as buildups made by forward/backward finite 
difference formul\ae. This approach is self-understanding
and we call it the standard approach. The result are basic schemes of
discretization. An illustration of numerical neighbourhoods of basic schemes
is given in Figure \ref{fig4.1}. These schemes are studied thoroughly in
\cite{SMMM}. The strict ellipticity of forms $a_n(\cdot,\cdot)$ is expressed always
in terms of the original pair of positive numbers $\Md, \Mg$. Generally,
the strict ellipticity of discretized forms follow from  the compartmental 
structure of $A_n$. These problems are analyzed in the next section.
\subsection*{Extended schemes}
Contrary to the standard approach in developing discretizations of $A_0(\msx)$
we have methods based on some general principle and which are not a priori
related to the forward/backward finite difference formulas. One of such methods
is described here. Principles of construction are given in \cite{LR3}.
To each $D_l$ we must associate elements $\msr(l) \in {\bbN}$ playing the role
analogous to $\msp$ for basic schemes. 

Again we assume the strict separation of sets $F(\pm)=\supp(\max\{\pm a_{12},
0\})$. The set $\EuScript{L}$ of Assumption \ref{ass4.1} is partitioned into 
the subsets $\EuScript{L}_{\mp}$, where $l \in \EuScript{L}_-$ means that 
$a_{ij} \leq 0$ on $D_l$ and $l \in \EuScript{L}_+$ means $a_{ij} \geq 0$ on $D_l$. 
Let us remind that the sets $D_l$ in present case are determined in terms of cubes
$C_{n(\varepsilon)}(\msr(l),\msx)$. There is always a room of arbitrariness in
a determination of these sets. The following maximal property removes
some of arbitrariness. There exist a $n(\varepsilon) \in {\bbN}$ such that
the sets $D_l, l \in \EuScript{L}_-$ have the following properties:
\begin{itemize}
 \item[a)] $a_{12} \leq 0$ on $\cup_{l \in \EuScript{L}_-}D_l$, 
 \item[b)] Each $D_l$ satisfies Assumption \ref{ass4.1},
 \item[c)] Each $C_{n(\varepsilon)}(\msr(l),\msx), l \in \EuScript{L}_-, \msx
 \in G_n$ on which $a_{12} = 0$ is contained in some of $D_l,l \in \EuScript{L}_-$.
\end{itemize}
Now we can define $D(-) = \cup_{l \in \EuScript{L}_-}D_l$ and $D(+)=
D(-)^c$. In our proceeding discussion we assume $n(\varepsilon)=1$. The
subgrids $G_n(-)=G_n \cap D(-)$ have the same properties as the corresponding
subrids $G_n(P,-), p_1=p_2=1$ in the subsection on basic schemes. Similarly, the
sets $G_n(+) = G_n\cap \overline{D(+)}$ coincide with $G_n(P,+)$ as well.
Now we define subgrids $G_n(l,-) = G_n(-) \cap D_l$ and conclude that
$G_n(l,-)$ form a partition of $G_n(-)$. However, the sets $G_n(l,+)= G_n(+) \cap D_l$
do not form a partition of $G_n(+)$ because some of grid-knots at $\partial 
D(+)$ may be outside of each $G_n(+,l)$. Therefore, we have to extend the sets
$D_l$ to wider sets $\tilde{D}_l$ such that $G_n(l,+)= G_n(+) \cap \tilde{D}_l$
form a partition of $G_n(+)$. The sets $\tilde{D}_l$ cannot be defined
uniquely. Here we demand the following properties. The sets must
be disjoint, and the closure of $int(\tilde{D}_l)$ must coincide with
$cls(D_l)$. Thus we have:
\[ G_n(l,-) \ = \ G_n(-) \cap D_l, \quad G_n(l,+) \ = \ G_n(+) \cap 
 \tilde{D}_l.\]
In accordance with our discussion about properties of sets $G_n(P,\pm)$ in
the subsection on basic schemes we finally conclude that $G_n(l,\pm)$ cover $G_n$
and some of them may have common grid-knots.

The forms $a_n(\cdot,\cdot)$ are expressed in terms of
$a_{ij}(\msx+\msx^{(\pm \pm)}(l,n)), \msx \in G_n$, where $\msx^{(\pm \pm)}
(l,n)$ are certain elements  in ${\bbR}^d$. Since $\msx^{(\pm \pm)}
(l,n)$ can take values outside of $cls(D_l)$ we are due to 
specify how to take values of $a_{ij}(\msx+\msx^{(\pm \pm)}(l,n))$ in 
such cases. The values must be taken in the set $D_l(\varepsilon)$ of Assumption 
\ref{ass4.1}. In this way we conclude that the entries of $A_n$ are calculated 
in terms of values of 
coefficients $a_{ij}$ at points which are not grid-knots. 

It is convenient to use a representation $a_n(v,u) = a_n^{(-)}(v,u) + 
a_n^{(+)}(v,u)$, where the forms $a_n^{(\mp)}(v,u)$ are related to the 
index sets $\EuScript{L}_{\mp}$ as previously. Let us define
\[\begin{array}{l}
 \mst^{(\pm +)}(\msr) \ = \ 
 \frac{1}{2}\,\big (\pm h\,r_i(l)\mse_i \, 
 +\,h\,r_j(l)\mse_j\big) \ \in S_n(\msr,{\bf 0}),\\
 \mst^{(\pm -)}(\msr) \ = \ 
 \frac{1}{2}\,\big (\pm h\,r_i(l)\mse_i \, 
 -\,h\,r_j(l)\mse_j\big) \ \in S_n(\msr,{\bf 0}).
 \end{array}\]
Obviously  $\msx^{(\alpha \beta)}(l,n)$ and $\mst^{(\alpha \beta)}(\msr)$
coincide for $\msp =\msr(l)$. The form $a_n^{(-)}(\cdot,\cdot)$ 
is defined by:
\begin{equation}\label{exp4.9}\begin{array}{l} \displaystyle
 a_n^{(-)}(v,u) \ = \ \sum_{l \in \EuScript{L}_-}\:\sum_{\msxd \in G_n(l)}\:\Big(
 \sum_{i = 1}^2\:
 a_{ii}(\msx+\mst^{(++)}({\bf 1}))\, \big(\der_i(h)v \big)(\msx) \big(
 \der_i(h)u \big)(\msx)
 \\  \displaystyle
 +\:\sum_{i \ne j}\:a_{ij}(\msx+\mst^{(++)}(\msr))\, \big(\der_i
 (r_i(l)h)v \big)(\msx)
 \big(\der_j(r_j(l)h)u \big)(\msx) 
 \\  \displaystyle 
 +\: \sum_{i \ne j}\:a_{ij}(\msx+\mst^{(++)}(\msr))\, \frac{r_i(l)
 h}{r_j(l)h}\:
 \Big[ \big(\der_i(h)v \big)(\msx) \big(\der_i(h)u \big)(\msx)\\
 \:-\:\big(\der_i(r_i(l)h)v \big)(\msx) \big(\der_i(r_i(l)h)u
 \big)(\msx)\Big] \ \Big). \end{array}
\end{equation}
We obtain $a_n^{(+)}(v,u)$ from $a_n^{(-)}(v,u)$ by replacing 
$\der_i(r_i(l)h)$ with $\der_i((-1)^{i-1}(r_i(l)h))$, $\EuScript{L}_-$ with
$\EuScript{L}_+$ and $\mst^{(++)}$ with $\mst^{(+-)}$. The forms $a_n^{(\mp)}(v,u)$
are not second degree polynomials of $\der_i(h)$ with
simple structure. Due to the compartmental structure of $A_n$ they can be
ultimately represented as forms depending on $\der_i(q_ih)$ with various
$q_i$. For the quantities $a_n^{(\mp)}(u,u)$ more comprehensible
expressions can be written down such as~(\ref{exp5.16}) and (\ref{exp5.17}).

Discretizations $A_n$ have a general expression:
\[ A_n \ = \ \sum_{l \in \EuScript{L}_-}\:A_n^{(-)}(l)\,{\bbJ}_{G_n(l,-)}\:+\:
 \sum_{l \in \EuScript{L}_+}\:A_n^{(+)}(l)\,{\bbJ}_{G_n(l,+)},\]
where
\begin{equation}\label{exp4.22}\begin{array}{l} \displaystyle
 A_n^{(-)}(l) 
 = \ -\:\sum_{i =1}^2\:V_i\,\EuScript{A}_n^{(+)}(i,i,1,1)\,U_i
 \:-\:\sum_{i \ne j}^2\: V_i(r_i(l))\, \EuScript{A}_n^{(+)}(i,j,r_i,r_j)
 \,U_j(r_j(l))
 \\ \displaystyle
 -\:\sum_{i \ne j}^2\:\frac{r_i(l) h}{r_j(l)h}\,
  \Big[V_i\, \EuScript{A}_n^{(+)}(i,j,1,1)\,U_i \:-\:
 V_i(r_i(l))\,\EuScript{A}_n^{(+)}(i,j,r_i,r_j)\,U_i(r_i(l)) \Big],\\ \\
 A_n^{(+)}(l) = \ -\:\sum_{i =1}^2\:V_i\,\Lambda_i(r_i)^T\,\EuScript{A}_n^{(-)}
 (i,i,1,1)\, \Lambda_i(r_i)\,U_i\\ \displaystyle
 -\:\sum_{i \ne j}^2\: V_i(r_i(l))\,\Lambda_i(r_i)^T \EuScript{A}_n^{(-)}
 (i,j,r_i,r_j) \Lambda_j(r_j)\,U_j(r_j(l))
 -\:\sum_{i \ne j}^2\:\frac{r_i(l) h}{r_j(l)h}\, 
 \\ \displaystyle
 \times \,  \Big[V_i\,\Lambda_i(r_i)^T\, \EuScript{A}_n^{(-)}(i,j,1,1)\,
 \Lambda_i(r_i)\,U_i \:-\:
 V_i(r_i(l))\,\EuScript{A}_n^{(-)}(i,j,r_i,r_j)\,U_i(r_i(l)) \Big], \end{array}
\end{equation}
and where $\Lambda_1(r_1) = I, \Lambda_2(r_2) =Z_n(-r_2(l),2)$ as in the previous case.

In order to write down the entries of $A_n$ we need the following abbreviations:
\[\begin{array}{l}
 \msw^{(\pm)}(l) \ = \ r_i(l)\mse_i \,\pm\,r_j(l)\mse_j \ \in I_n,\\
 a_{ij}^{(\alpha \beta)}(\msr) \ = \ a_{ij}( \msx + \mst^{(\alpha \beta)}(\msr)),
 \quad \alpha, \beta \in \{+,-\},\\
 \hat{a}_{12}^{(-+)}(\msr) \ = \ a_{12}(\msx+\mst^{(++)}(\msr)-
 h\mse_1),\\
 \hat{a}_{12}^{(+-)}(\msr) \ = \ a_{12}(\msx+\mst^{(++)}(\msr)-
 h\mse_2),\\
 \hat{a}_{ii}^{(++)}(\msr) \ = \ a_{ii}^{(++)}(\msr),\quad
 \hat{a}_{ii}^{(--)}(\msr) \ = \ a_{ii}^{(--)}(\msr).
\end{array}\]
Then we have the following nontrivial off-diagonal entries of $A_n$:
\begin{equation}\label{exp4.23}\begin{array}{ll}\displaystyle
 \big(A_n\big)_{\mskd \mskd\pm \msed_1}& \displaystyle
 = \  -\,\frac{1}{h^2}
 \: \left\{\begin{array}{lll} \displaystyle
 a_{11}^{(\pm +)}({\bf 1})\:-\: \frac{r_1(l)}{r_2(l)}\,
 \big|\hat{a}_{12}^{(\pm +)}(\msr)\big| &{\rm for}& a_{12} \leq 0, \\ 
 \displaystyle
 a_{11}^{(\pm -)}({\bf 1})\:-\: \frac{r_1(l)}{r_2(l)}\,
 \big|\hat{a}_{12}^{(\pm -)}(\msr)\big| &{\rm for}& a_{12} \geq 0,  \end{array}
 \right . \\ \displaystyle
 \big(A_n\big)_{\mskd \mskd\pm \msed_2} & \displaystyle
 = \  -\,\frac{1}{h^2}
 \: \Big[ a_{22}^{(+\pm)}({\bf 1})\:-\: \frac{r_2(l)}{r_1(l)}\,
 \big|\hat{a}_{12}^{(+\pm)}(\msr)\big|\Big].  \end{array}
\end{equation}
The entries in the plane spanned by $\mse_1, \mse_2$ have the structure:
\begin{equation}\label{exp4.24}\begin{array}{l}\displaystyle
 \big(A_n\big)_{\mskd \mskd \pm \mswd^{(-)}(l)} \ = \  
 -\,\frac{1}{h^2r_1(l)r_2(l)}\:
 \big|a_{12}^{(\pm \mp)}(\msr)\big| \quad {\rm for}\quad a_{12} \leq 0 \ {\rm on} \ D_l,
 \\ \displaystyle
 \big(A_n\big)_{\mskd \mskd \pm \mswd^{(+)}(l)} \ = \  
 -\,\frac{1}{h^2r_1(l)r_2(l)}\:
 \big|a_{12}^{(\pm \pm)}(\msr)\big| \quad {\rm for} \quad a_{12} \geq 0 \ {\rm on} \ D_l,
 \end{array}
\end{equation}
where $+\msw^{(\pm)}$ is associated with the upper and $-\msw^{(\pm)}$ with the
lower indices of $a^{(\alpha \beta)}$, respectively. 

\begin{figure}
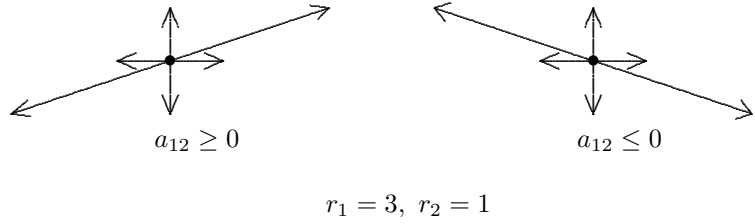

\caption{Numerical neighbourhoods at internal grid-knots of $D_l$}
\label{fig4.2}
\beginpicture
\setcoordinatesystem units <1.0pt,1.0pt>
\setplotarea x from -40 to 300, y from 0 to 100
\setlinear
\plot 40 60, 80 60 /
\plot 60 40, 60 80 /
\plot 0 40, 120 80 /
\plot 200 60, 240 60 /
\plot 220 40, 220 80 /
\plot 160 80, 280 40 /
\plot 47 57, 40 60, 47 63 /
\plot 73 57, 80 60, 73 63 /
\plot 7 39, 0 40, 6 45 /
\plot 116 76, 120 80, 113 81 /
\plot 57 47, 60 40, 63 47 /
\plot 57 73, 60 80, 63 73 /
\put{$\bullet$} [c] at 60 60
\put{$\bullet$} [c] at 220 60
\plot 207 57, 200 60, 207 63 /
\plot 233 57, 240 60, 233 63 /
\plot 217 47, 220 40, 223 47 /
\plot 217 73, 220 80, 223 73 /
\plot 166 82, 160 80, 164 75 /
\plot 274 38, 280 40, 276 45 /
\put{$r_1 = 3, \ r_2 =1$} [c] at 150 5
\put{$a_{12} \geq 0$} [c] at 70 30
\put{$a_{12} \leq 0$} [c] at 230 30
\endpicture
\end{figure}

Discretizations $A_n$ defined by (\ref{exp4.23}), (\ref{exp4.24}) are called 
extended schemes. The numerical neighbourhoods are illustrated in
Figure~\ref{fig4.2}.

In (\ref{exp4.22}) we have 4 sums with respect to the indices $i,j$. The
first and second sums have expressions similar to Expressions (\ref{exp4.5}).
They contribute to Expressions (\ref{exp4.24}) and to a part of entries in
(\ref{exp4.23}). Unfortunately there appear non-trivial entries $(A_n)_{
\mskd \msld}, \msl = \msk \pm r_i\mse_i$. These entries must be canceled
by contributions from the third and fourth sums. These two latter sums lack
the structure similar to (\ref{exp4.5}) since the sum includes the terms
$a_{12}\der_1v \der_1u$ and $a_{12}\der_2v \der_2u$. So
the net result of all four sums are entries (\ref{exp4.24}).
Let $J \subset I_n(l)$ consists of grid-knots in some of $\msx \in G_n(l,\pm)$ 
for which the numerical neighbourhoods ${\cal N}(\msx)$ have only internal
grid-knots, ${\cal N}(\msx) \subset int(G_n(l,\pm))$. The associated
diagonal submatrix $(A_n)_{J J}$ is symmetric. Diagonal submatrices for
which ${\cal N}(\msx)$ contain boundary grid-knots of $G_n(l,\pm)$ may
lack the symmetry. If the quantities $\hat{a}_{ij}^{(\alpha \beta)}$ in 
(\ref{exp4.23}) are replaced with $a_{ij}^{(\alpha \beta)}$, the symmetry
of $A_n$ is lost altogether, although the convergence is still
preserved. However, the quantities $a_{ij}^{(\alpha \beta)}$ must
not be replaced with $a_{ij}(h\msk)$ since the resulting $(A_n)_{\mskd \msld}$
would be discretizations of $-\sum_{ij}a_{ij}\partial_i \partial_j$.

\subsection{Discretization with compartmental structure}
Now we can describe general structure of constructed discretizations of 
$A_0(\msx)$. From the definition of bilinear forms $a_n(u,v) = \lev {\bf v}|A_n
{\bf u}\des$ the following property is obvious: If ${\bf v} = {\bf 1}$ and
${\bf u}$ are with compact supports on $G_n$ there must be $a_n(1,u)=a_n(u,1)=0$,
implying that the row sums and column sums of $A_n$ have zero values. Hence,
if the off-diagonal entries of $A_n$ are non-positive then the matrices 
$A_n$ are simultaneously $M$-matrices and have the compartmental structure.
We consider here only the extended scheme.

It is convenient to utilize the quantities:
\begin{equation}\label{exp4.7}
 \omega_n(a_{ii},\msx) \ = \ \frac{1}{d-1}\:\sum_{s \ne i}^d \: 
 a_{ii}(h\msk+h\msm_{ii}(l,s)),
\end{equation}
where $\msm_{ii}(l,s)$ are defined by the rules of construction of
extended schemes.

\begin{procedure}\label{dsp4.1} Let Assumption \ref{ass4.1} be valid and
matrices $A_n$ on $G_n$ be constructed by the rule (\ref{exp4.14}). Then
their entries have the following properties:
\begin{enumerate}
 \item Entries of $(A_n)_{\mskd \msld}, \msk, \msl \in I_n$, $\msx = h \msk$ are
linear combinations of $a_{ij}(\msx_{ij}(n,\msx,l))$ where $\msx_{ij}(
n,\msx,l) = h\msk + h\msm_{ij}(l,s)$, $\msm_{ij}(l,s)$ are $n$-independent 
elements of ${\bbR}^d$ and $l \in \EuScript{L}, s \in \{1,2,\ldots,d\}$.
 \item For each grid-knot $\msx = h \msk$: $(A_n)_{\mskd \mskd} = -
\sum_{\msld}  (A_n)_{\mskd \msld}$.
 \item For each $\msx = h \msk \in cls(D_l)$ entries on the "cross branches" $\msx 
\pm h\mse_i$, i.e. $(A_n)_{\mskd \pm h\msed_i}$, are 
defined by:
\[\begin{array}{lll} \displaystyle
 \big(A_n \big)_{\msxd \msxd\pm h\msed_i} &=& \displaystyle
 -\,\frac{1}{h^2}  \,\Big [ \omega_n(a_{ii},\msx) \\ 
 \displaystyle && \displaystyle
 -\,\sum_{m \ne i}\:\frac{r_i(l)}{r_m(l)}\, 
 |a_{im}(\msx_{im}(n,\msx,l))|\Big]. \end{array}\]
 \item Entries of $A_n$ which are not on the "cross branches" are defined by
using elements $\msz_{ij}(l) = r_i(l)\mse_i - r_j(l)\mse_j \in I_n$ or elements
$\msz_{ij}(l) = r_i(l)\mse_i + r_j(l)\mse_j \in I_n$:
\[ \big(A_n\big)_{\mskd \mskd \pm \mszd_{ij}(l)} \ = \ -\:\frac{1}{
 r_i(l)  r_j(l)} \: |a_{ij}(\msx_{ij}(n,\msx,l))|.\]
\end{enumerate}
\end{procedure}

Some peculiar features regarding the structure of sets ${\cal N}(\msx), \msx = h\msk \in 
G_n \cap D_l$, must be pointed out. If $a_{ij}, i \ne j$ is not changing the sign
in a neighbourhood of $\msx$ then the minimal
number of elements in ${\cal N}(\msx)$ is $1+d+d^2$. In this case the set 
${\cal N}(\msx)$ consists of its center, $2d$-grid-knots on the $d$-dimensional
cross $\{\pm \mse_i: i = 1,2,\ldots,d\}$ and 2 grid-knots in each two-dimensional
plane. Generally, the number of grid-knots in two-dimensional plane may be
larger than 2. Here we consider only the case of two grid-knots at most in 
the two-dimensional planes. 
This demand has the following implication on the construction of discretizations
$A_n^{(rs)}$. Let the pairs $\mse_r, \mse_s$ and $\mse_s, \mse_t$ define
two-dimensional planes and let $A_n^{(rs)}, A_n^{(st)}$ be the corresponding
discretizations which are constructed using the parameters 
$\msr^{(rs)}, \msr^{(st)}$. Then there must hold $(\msr^{(rs)})_s =
(\msr^{(st)})_s$. In such case the off-diagonal entries of $A_n$ 
have the structure as described in 4. of Discretization procedure \ref{dsp4.1}.
The described structure of sets ${\cal N}(\msx)$ is valid for all $\msx$
because the functions $a_{ij}, i \ne j$ do not change sign on $S_n(\msp,\msx)$. 
This is an important consequence of the strict separation of sets $\max\{a_{12},0\}$. 

Obviously that all the constructed forms $a_n$ of this section are discretizations
of the form (\ref{exp2.4}).
One can be easily convinced that the terminology "a discretization of the original
form (\ref{exp2.4})" is not artificial. At the present level of analysis it is easy to
check $a(v,u) = \lim_n h^d a_n(v,u)$ for any pair $v,u \in C_0^{(1)}({\bbR}^d)$.

\begin{theorem}\label{Th4.1} Let Assumption \ref{ass4.1} be valid. There exist
discretizations $A_n$ which are constructed by the rules of Discretization procedure
\ref{dsp4.1}, such that $A_n$ have the compartmental structure.
\end{theorem}

{\Proof} For each $D_l$ we have to choose the parameters $r_i(l)$ of the properties
3. and 4. of Discretization procedure \ref{dsp4.1} so that the condition 3. of
Assumption \ref{ass4.1} is valid. The uniform continuity of coefficients on $D_l$ and
the condition 3. of Assumption \ref{ass4.1} ensure the compartmental structure of
matrices $A_n$ as demonstrated in \cite{LR3}. {\QED}

\section{Convergence in $W_2^1$-spaces}\label{sec5}
Discretizations of the original variational problem (\ref{exp2.9}) or (\ref{exp2.10}) 
are defined in terms of a sequence of bilinear forms $a_n(\cdot,\cdot)$ on 
$E_n(R,{\bbR}^d) \times E_n(R,{\bbR}^d)$ and a sequence of linear functionals 
$\lev \cdot |\msmu_n\des_R$ on $E_n(R,{\bbR}^d)$. The
associated discretized variational problems are defined by equalities
(\ref{exp4.1}). The discretized 
variational equalities (\ref{exp4.1}) can be rewritten in an equivalent manner:
\begin{equation}\label{exp5.3}
 (\lambda I \,+\,A_n)\,{\bf u}_n \ = \ \msmu_n, 
\end{equation}
where in the case of Problem (\ref{exp2.10}) ${\bf u}_n, \msmu_n$ are grid-functions 
on $G_n$ and in the case of Problem (\ref{exp2.9}) they are grid-functions on 
$G_n(R,D)$ or $G_n(R,D)$.

Problems (\ref{exp2.9}) or (\ref{exp2.10}) are solved numerically in two steps.
In the first step we construct grid-functions ${\bf u}_n$ on $G_n(R)$ or
$G_n(R,D)$ according to~(\ref{exp5.3}). The obtained grid-functions represent
the solution at grid-knots, therefore, we call them {\em grid-solutions}. 
The grid-solutions are imbedded
into the spaces $E_n(R,{\bbR}^d)$ or $E_n(R,D)$ by~(\ref{exp3.10}), and the
convergence $u(n) \to u$ must be proved in some Banach spaces. Therefore we call
functions $u(n)$ {\em approximate solutions}.

Though the functions $u(n) = \Phi(R) {\bf u}_n$ are called approximate
solutions, this terminology has to be justified after a convergence
analysis. The convergence proofs are based on some properties of the forms
$a_n(\cdot,\cdot)$ and linear functionals $\lev \cdot |\msmu_n\des_R$ to be 
described in details later in this section. Most of the analysis in this section is
related to Problem (\ref{exp2.10}) on ${\bbR}^d$. The obtained results can be easily 
applied to Problem (\ref{exp2.9}) on a bounded domain. This is carried out at the
end of section.
\subsection{Consistency}
Certain number of notions important for the convergence of approximate solutions
is formulated in terms of sequences of functions with a particular structure:
\begin{equation}\label{exp5.7}\begin{array}{c} \displaystyle
 \mathfrak{V} \ = \  \{v(n)\::\: n \in {\bbN}\} \subset \cup_n E_n(R,{\bbR}^d),\\
 \mathfrak{U} \ = \  \{u(n)\::\: n \in {\bbN}\} \subset \cup_n E_n(R,{\bbR}^d).
 \end{array}
\end{equation}
\begin{definition}[Consistency]\label{def4.3} The forms $a_n(\cdot,\cdot)$
on $E_n(R,{\bbR}^d) \times E_n(R,{\bbR}^d)$ are consistent with the form (\ref{exp2.4}) 
if
\[ a(v,u) \ = \ \lim_n h^d\,a_n(v(n),u(n)) \]
is valid for any pair $\mathfrak{V}, \mathfrak{U}$ of (\ref{exp5.7}) such
that $\mathfrak{V}$ converges weakly in $W_2^1({\bbR}^d)$ to $v$, and
$\mathfrak{U}$  converges strongly in $W_2^1({\bbR}^d)$ to $u$.
\end{definition}
\begin{proposition}\label{Prop5.1} Let a sequence of matrices $\{A_n:n \in {\bbN}\}$
be constructed by basic schemes or extended schemes and let $\{a_n(\cdot,\cdot):
n \in {\bbN}\}$ be the corresponding sequence of discretizations of (\ref{exp2.4}).
Then the forms $a_n$ are consistent with the form (\ref{exp2.4}).
\end{proposition}
This important result is proved by a lemma which is formulated bellow.
Let $\mathfrak{V}$ and $\mathfrak{U}$ be defined by (\ref{exp5.7}) and
converge weakly and strongly in $W_2^1({\bbR}^d)$ to $v, u \in W_2^1({\bbR}^d)$,
respectively. Expressions $\omega_n(a_{ij},\msx)$ for $i=j$ are defined by
(\ref{exp4.7}). In this proof we extend this definition to the case $i \ne j$ 
and define $\omega_n(a_{ij},\msx) = a_{ij}(\msx_{ij}(n,\msx,l))$,
where $\msx_{ij}(n,\msx,l) = \msx + h\:\msm_{ij}(l,n,s)$ are constructed 
by the rules of basic and extended schemes. For each 
$l \in \EuScript{L}$ and each pair $i,j \in \{1,2,\ldots,d\}$ there holds:
\[\begin{array}{c} \displaystyle
 \lim_n \:\int_{D_l(\varepsilon)} \: a_{ij}(\msx)\,(\partial_iv(n))(\msx)\, 
 (\partial_ju(n))(\msx)\, d\msx \\ \displaystyle
  = \ \lim_n\:\sum_{\msxd \in G_n(P) \cap D_l(\varepsilon)}\:\omega_n(
 a_{ij},\msx)\,\int_{C_n(\mspd(l),\msxd)}\:
 (\partial_iv(n))(\msx)\, (\partial_ju(n))(\msx)\,d\msx,\end{array}\]
Let us point out that the written identity is
unchanged if we replace $\omega_n(a_{ij},\msx)$ with $a_{ij}(\msx)$.
Let us consider a bilinear form on $W_2^1({\bbR}^d) \times W_2^1({\bbR}^d)$ 
defined by
\begin{equation}\label{exp5.8}
 \theta_n(l,v,u) \ = \ \sum_{\msxd \in G_n(R)\cap D_l} \:p(
 \msw(n,\msx))\,\int_{C_n(\msrd,\msxd)}\:\partial_i v(\msy)\,\partial_j u(\msy)
 d \msy,
\end{equation}
where $p$ is a uniformly continuous function on $D_l$ and
$\msw(n,\cdot)$ is a transformation of ${\bbR}^d$ such that $\msw(n,\msx) \in 
C_n(\msr,\msx)$. Then
\begin{equation}\label{exp5.9}
 \lim_n \,\theta_n(l,v(n),u(n)) \ = \ \int_{D_l} \:p(\msy)\,
 \partial_i v(\msy)\,\partial_j u(\msy) d \msy 
\end{equation}
for any pair $\mathfrak{V}, \mathfrak{U}$ of (\ref{exp5.7}), converging in 
$W_2^1({\bbR}^d)$ weakly to $v$ and strongly to $u$, respectively. 

The object of next analysis is the bilinear functional on $E_n({R,\bbR}^d) \times
E_n({R,\bbR}^d)$ defined by:
\begin{equation}\label{exp5.10}
 \gamma_n(l,v,u) \ = \ \ h^d\,\vol(R)\,\sum_{\msxd \in G_n(R)\cap D_l} \:p(
 \msw(n,\msx))\, \der_i(r_ih) v(h\msk)\, \der_j(r_jh) u(h\msk).
\end{equation}
In particular:
\[ \gamma_n(l,v(n),u(n)) \ = \ \ h^d\,\vol(R)\,\sum_{\msxd \in G_n(R)\cap D_l} \:p(
 \msw(n,\msx))\,v_{\mskd i}\,u_{\mskd j},\]
where $v_{\mskd i}=\der_i(r_ih) v(h\msk), u_{\mskd j}=\der_j(r_jh) u(h\msk)$.

\begin{lemma}\label{lem5.4} Let the sequence $\mathfrak{V}$ converge weakly in
$W_2^1({\bbR}^d)$ to $v$ and $\mathfrak{U}$ converge strongly to $u$ in 
$W_2^1({\bbR}^d)$. Then
\[ \lim_n \:\gamma_n(l,v(n),u(n)) \ = \ \lim_n \:\theta_n(l,v(n),u(n)) \ = \ 
\int_{D_l}\: p(\msx)\,\partial_i v( \msx)\, \partial_j u(\msx)\:d\msx.\]
\end{lemma}

{\Proof} We have to analyze $\theta_n(l,v(n),u(n))$ as $n \to \infty$. 
To simplify notation we assume $i, j \in \{1,2\}$. 
The indices $\msk = (k_1,k_2,\ldots,k_s)$ are denoted shortly as $\msk = (k_1,
\msk')$ and $\msk = (k_1,k_2,\msk'')$, $k_1,k_2 \in {\bbZ}$. The corresponding
$r_1, r_2$ are denoted by $r,s$, respectively. After inserting expressions for
$v(R,n), u(R,n)$ into (\ref{exp5.8}) and carrying out a straightforward calculation
we get expressions:
\[ \theta_n(l,v(n),u(n)) \ = \ \sum_{\msxd \in G_n(R) \cap D_l}\:
 p(\msw(n,\msx))\,\rho_n(h\msk,v(n),u(n)),\]
where $\rho_n(h\msk,v(n),u(n))$ denotes an integral which for the case of $i = j =1$
has the form:
\[\begin{array}{c}\displaystyle
 \rho_n(h\msk,v(n),u(n)) \ = \ 
 \int_{C_n(\msrd,h\mskd)}\: \partial_1 v(R,n)(\msy)\,\partial_1 u(R,n)(\msy)
 \,d \msy \ = \ \sum_{\mskd' \msld'} \:(\psi_{\mskd'}|\psi_{\msld'})\\ \displaystyle
 \times \:\left [  \sum_{l=k,k+r}\,v_{k \mskd'} u_{l \msld'}\:
 (\partial_1 \psi_k |\partial_1 \psi_l) \:+\:
 \sum_{l=k,k+r}\,v_{k+r \mskd'} u_{l \msld'}\: (\partial_1 \psi_{k+r} |\partial_1 \psi_l)
 \right] \\  \displaystyle
 = \ rh\:\sum_{\mskd' \msld'} \:(\psi_{\mskd'}|\psi_{\msld'})\:
 \big(\der_1(rh)v\big)(h \msk)\, \big(\der_1(rh)u\big)(h \msk).  \end{array} \]
Now we use quantities $s_{\mskd \msld} = (\psi_{\mskd}|\psi_{\msld})\Vert \psi_{\mskd}
\Vert_1^{-1}$ and bring into mind their properties $s_{\mskd \msld} \geq 0, 
\sum_{\msld} s_{\mskd \msld} = 1$. 
\[ \theta_n(l,v(n),u(n))  =  h^d\,\vol(R)\:\sum_{k,\mskd',\msld'\in I_n(R),
\msxd \in D_l}\: 
 p(\msw(n,\msx)) \: s_{\mskd' \msld'}\:  v_{(k\mskd') i}\,u_{(k \msld') i}.\]
Analogously we get for $i=1, j=2$:
\[\begin{array}{c} 
 \theta_n(l,v(n),u(n)) \ = \ h^{d-2}\,\vol(R'')\!\!\!\!\!
 \sum\limits_{k,l,\mskd'',\msld''\in I_n(R), \msxd \in D_l}\: 
 s_{\mskd'' \msld''}\: p(\msw(n,h \msk))\: \int\limits_{J(k,r) \times
 J(l,s)}\!\!\!\!\!\!\!dz_1 dz_2 \\ \displaystyle \times \:
 \Big[v_{(k,l,\mskd'') i}\,u_{(k,l, \msld'') j}\,\psi_k(z_1)\psi_l(z_2) \,+\,
 v_{(k,l,\mskd'') i}\,u_{(k+r,l, \msld'') j}\,\psi_{k+r}(z_1)\psi_l(z_2) \\
 \displaystyle  +\:
 v_{(k,l+s,\mskd'') i}\,u_{(k,l, \msld'') j}\,\psi_k(z_1)\psi_{l+s}(z_2) \,+\,
 v_{(k,l+s,\mskd'') i}\,u_{(k+r,l, \msld'') j}\,\psi_{k+r}(z_1)\psi_{l+s}(z_2)
 \Big], \end{array}\]
where $J(k,r) = [hk,h(k+r)]$ and $R''$ stands for the parameter set $\{\msr_0,r_3,
r_4,\ldots,r_d\}$. Upon integration over $z_1, z_2$ we get
\[\begin{array}{c} \displaystyle
 \theta_n(l,v(n),u(n)) = \frac{1}{4}\:h^{d-2}\,\vol(R'')\:
 \sum_{k,l,\mskd'',\msld''\in I_n(R), \msxd \in D_l}\: 
 s_{\mskd'' \msld''} p(\msw(n,h \msk))\\ \displaystyle
 \Big[v_{(k,l,\mskd'') i}\,u_{(k,l, \msld'') j}+
 v_{(k,l,\mskd'') i}u_{(k+r,l, \msld'') j}+
 v_{(k,l+s,\mskd'') i}u_{(k,l, \msld'') j}+
 v_{(k,l+s,\mskd'') i}u_{(k+r,l, \msld'') j} \Big]. \end{array}\]
We finish the proof for the case $i = j =1$ since the case $i \ne j$ can be treated
analogously.
The quantity $\gamma_n$ would be equal to $\theta_n$ if $s_{\mskd \msld}$ were absent 
and the double sum were replaced with the single sum over indices $\msk$.
Therefore we are due to estimate their difference:
\begin{equation}\label{exp5.11}
 \gamma_n(l,v(n),u(n))\,-\,\theta_n(l,v(n),u(n)) \ = \
 h^d\,\vol(R)\,\sum_{\mskd \msld} \:p(\msw(n,h\msk))\,v_{\mskd i}\,\delta_{kl}\,
 s_{\mskd' \msld'} \,\big(u_{\msld j}-u_{\mskd j}\big).
\end{equation}
To estimate the right hand side we need $\overline{p} = \sup p$:
\[\begin{array}{c} \displaystyle
 \Big|\gamma_n(l,v(n),u(n))\,-\,\theta_n(l,v(n),u(n))\Big| \ \leq \\ \displaystyle
 \overline{p}\:h^{d/2}\lnorm U_i(r_i) {\bf v}_n \lnorm_{R2} \:h^{d/2}\,
 \max \{\lnorm \big(Z(r_j,j)-I\big) U_j(r_j){\bf u}_{n}\lnorm_{R2} \,:\,
 j = 1,2,\ldots,d\}. \end{array}\]
By Theorem \ref{thn3.1} we have
\[\begin{array}{l}
 h^d\,\lnorm U_i(r_i){\bf v}_{n}\lnorm_{R2}^2 \ \leq \ (1-\sigma^2)^{-1}\,
 \Vert \,\partial_i\,v(R,n)\Vert_2^2,\\ 
 h^d\,\lnorm (Z(r_j,j)-I)\,U_j(r_j){\bf u}_{n}\lnorm_{R2}^2 \ \leq \
 (1-\sigma^2)^{-1}\, \Vert (Z(\msw)-I)\,\partial_j\,u(R,n)\Vert_2^2,\end{array}\]
where $\msw = hr_j \mse_j$. Due to the strong convergence of $\mathfrak{U}$ we have 
\[\begin{array}{c} \displaystyle
 h^{d/2}\,\lnorm \big(Z(r_j,j)-I\big) U_j(r_j){\bf u}_{n}\lnorm_{R2} 
 \ \leq \ (1-\sigma^2)^{-1/2}\:\Vert \big(Z(\msw)-I\big) \,\partial_j\,u(n,R)\Vert_2 \\
 \displaystyle
 \leq \ (1-\sigma^2)^{-1/2}\,\Big [\Vert (Z(\msw)-I) \,\partial_j\,u \Vert_2 \,+\,2\,
 \Vert \,\partial_j\,u - \,\partial_j\,u(R,n)\Vert_2 \Big] \ \to \ 0, \end{array}\]
so that $\lim_n \gamma_n(l,v(n),u(n)) = \lim_n \theta_n(l,v(n),u(n)) = 
\int p \partial_iv \partial_j u$. {\QED}

Now a proof of Proposition \ref{Prop5.1} follows from the inequalities
\[ \lim_n \: \sum_{l \in \EuScript{L}'}\: \gamma_n(l,v(n),u(n)) \ = \
 \int_{\cup \{D_l:l \in \EuScript{L}'\}}\: a_{ij}(\msx)\,  \partial_iv(\msx)
 \partial_ju(\msx)\,d \msx. \]
in which $p(\msw(n,h\msk))$ of (\ref{exp5.10}) is replaced with $a_{ij}(n,
\msx,l)$.

\subsection{Strict ellipticity of discretized forms}
A discrete form $a_n(\cdot,\cdot)$ on $l_0(G_n(R)) \times l_0(G_n(R))$ is said to be 
strictly elliptic \cite{Yo} if there exist two positive numbers $\Md(a_n), \Mg(a_n)$ 
such that
\[ \Md(a_n)\,\sum_{i=1}^d\:\lnorm U_i(r_i) {\bf u} \lnorm_{R2}^2 \ \leq \ 
 a_n(u,u) \ \leq \  \Mg(a_n)\,\sum_{i=1}^d\:\lnorm U_i(r_i) {\bf u} \lnorm_{R2}^2.\]
For a sequence of discrete forms $a_n(\cdot,\cdot)$ we need a stronger result.
The strict ellipticity must be uniform with respect to $n$ and $\vol(R)=\prod
r_i$ different values of $\msr_0$ in the parameter set $R=(\msr_0,\msr)$.

\begin{definition}\label{def5.1} Let $A_n$ be discretized by the rules of 
Discretization procedure \ref{dsp4.1}. Discrete forms $a_n(\cdot,\cdot)$ on 
$l_0(G_n(R)) \times l_0(G_n(R))$ are said to be strictly elliptic uniformly 
with respect to $n \in {\bbN}$ if there exist positive numbers $\Md \leq \Mg$
such that
\begin{equation}\label{exp4.4}
 \Md \,\sum_{i=1}^d\:\lnorm U_i(r_i) {\bf u} \lnorm_{R2}^2 \ \leq \
 a_n({\bf u}_n,{\bf u}_n)_R \ \leq \
 \Mg \,\sum_{i=1}^d\:\lnorm U_i(r_i) {\bf u} \lnorm_{R2}^2.
\end{equation}
for all $n \in {\bbN}$ and all $0 \leq (\msr_0)_i < r_i, i =1,2,\ldots,d$.
\end{definition}

\begin{proposition}\label{Prop4.2} Let the discretizations $A_n$ of
$A_0(\msx) = -\sum \partial_i a_{ij}(\msx)\partial_j$ be constructed by the rules
of Discretization procedure \ref{dsp4.1}. If $A_n$ have the compartmental 
structure then the discrete forms ${\bf v}, {\bf u} \mapsto \lev {\bf v}|A_n
{\bf u}\des_R$ are strictly elliptic on $l_0(G_n(R)) \times l_0(G_n(R))$ 
uniformly with respect to $n \in {\bbN}$.
\end{proposition}

{\Proof} First we consider the basic scheme for a two-dimensional grid. From 
the structure of bilinear forms (\ref{exp4.5}) we have
\begin{equation}\label{exp5.15}
 a_n(u,u) \ \geq \ \Md\,\lnorm U_1(p_1){\bf u}_n\lnorm_2^2 \,+\, 
 \Md\,\Big(\lnorm K(-)U_2(-p_2){\bf u}_n\lnorm_2^2+
 \lnorm K(+)\Lambda_2 U_2(-p_2){\bf u}_n\lnorm_2^2\Big),
\end{equation}
where $K(\mp)={\bbJ}_{G_n(P,\mp)}$ are projectors. Let us consider a pair 
$\msx, \msy=\msx-hp_2\mse_2$ which is involved in the definition
of operator $U_2(-p_2)$. The indices of $\msx, \msy$ are $\msk, \msl$, 
respectively. If $\msx, \msy \in G_n(P,+)$ then $(U_2(-p_2){\bf u}_n)_{\mskd}$ 
can be replaced with $(U_2(p_2){\bf u}_n)_{\msld}$. If $\msy \notin 
G_n(P,+)$ then $\msy \in G_n(P,-)$. In this case the term $(U_2(-p_2){\bf u}_n
)_{\msld}$ can be omitted from the sum in (\ref{exp5.15}) because this term 
is already contained in the corresponding sum of $U_1(p_1){\bf u}_n$. Hence, the 
right hand side of (\ref{exp5.15}) can be estimated from bellow by $\Md \sum_{i=1}^2
\lnorm U_i(p_i){\bf u}_n\lnorm_2^2$. In this way the left hand side of inequality
(\ref{exp4.4}) is proved. The right hand side of this double inequality 
follows by choosing the double value of $\Mg$.
In the case of $d > 2$ we use the construction (\ref{exp4.14}) 
and get the same lower and upper bounds. Let us point out that the compartmental 
structure is not used in this step of proof.

Let us now consider a two-dimensional problem with an extended scheme. In the
present step, the compartmental structure is utilized in an essential
way. The form $a_n^{(-)}(u,u)$ of (\ref{exp4.9}) can be rewritten as:
\begin{equation}\label{exp5.16}\begin{array}{l}
 a_n^{(-)}(u,u) \ = \\ \displaystyle
 \sum_{l\in {\cal L}_-}\,\sum_{\mskd \in G_n(l,-)}
 \sum_{i \ne j}\,\Big( a_{ii}(\msx+\mst^{(++)}({\bf 1}))+
 \frac{r_i(l)}{r_j(l)}a_{12}(\msx+\mst^{(++)}(\msr))\Big)
 (\der_i(h)u)(\msx)^2 \ -\\ \displaystyle
 \sum_{l\in {\cal L}_-}\,\sum_{\mskd \in G_n(l,-)}\!\!\!\!\!
 a_{12}(\msx+\mst^{(++)}(\msr))\,\left(
 \sqrt{\frac{r_1(l)}{r_2(l)}}(\der_1(r_1h)u)(\msx)\!-\!
 \sqrt{\frac{r_2(l)}{r_1(l)}}(\der_2(r_2h)u)(\msx)\right)^2.\end{array}
\end{equation}
Let us remind that $a_{12} \leq 0$ on the set $D(-)$.
Due to the compartmental structure the first term is positive definite.
The second term is positive semidefinite and can be
disregarded in the next step of estimation from bellow. The result
is:
\[ a_n^{(-)}(u,u) \ \geq \ \omega(a)\:
 \sum_{l\in {\cal L}_-}\,\sum_{\mskd \in G_n(l,-)}
 \sum_{i}\,\big(U_i{\bf u}_n\big)_{\mskd}^d.\]
where $\omega(a)$ is the positive number specified in Assumption \ref{ass4.1}.
For the form $a_n^{(+)}(u,u)$ we have an analogous inequality involving
the summation over all the indices $\msk \in G_n(+)$.
\begin{equation}\label{exp5.17}\begin{array}{l}
 a_n^{(+)}(u,u) = \\ \displaystyle
 \sum_{l\in {\cal L}_+}\sum_{\mskd \in G_n(l,+)}
 \sum_{i \ne j}\Big( a_{ii}(\msx+\mst^{(+-)}({\bf 1}))\!-\!
 \frac{r_i(l)}{r_j(l)}a_{12}(\msx+\mst^{(+-)}(\msr))\Big)
 (\der_i((-1)^{i-1}h)u)(\msx)^2 \!+\! \\ \displaystyle
 \sum_{l\in {\cal L}_+}\sum_{\mskd \in G_n(l,+)}\!\!\!\!\!
 a_{12}(\msx+\mst^{(+-)}(\msr))\left(
 \sqrt{\frac{r_1(l)}{r_2(l)}}\,(\der_1(r_1h)u)(\msx)-
 \sqrt{\frac{r_2(l)}{r_1(l)}}(\der_2(-r_2h)u)(\msx)\right)^2\!\!\!.\end{array}
\end{equation}
The lower bound follows in the same way as for $a_n^{(-)}$. {\QED}
\subsection{$W_2^1$-convergence}
We have shown how the form $\lambda (v|u)+a(v,u)$ is discretized by forms 
$h^d \lambda  \lev {\bf v}_n|{\bf u}_n\des_R+\lev {\bf v}_n|A_n{\bf u}_n\des_R$. 
In order to solve discretized problem~(\ref{exp5.3}), we
have to describe a discretization of the linear function $v \to \lev v|\mu\des$
by $h^d\lev {\bf v}_n|\msmu_n\des_R$, where $\msmu_n \in l(G_n(R))$.
First we must demonstrate the existence of $\msmu_n$ such that
$h^d\lev {\bf v}_n|\msmu_n\des_R \to \lev v|\mu\des$. Discretizations of
$\mu$ are defined by~(\ref{exp3.7}) so that (\ref{exp3.9}) is valid.
\begin{lemma}\label{lem5.2} Let $\mu$ be a continuous linear functional on 
$W_2^1({\bbR}^d)$. There exists discretizations $\msmu_n(R) \in l(G_n(R))$ such 
that 
\begin{equation}\label{exp5.5}
 \lev u(n) \,|\,\mu\des \ = \ h^d\:\lev {\bf u}_n\,|\,\msmu_n\des_R \ = \
 h^d\:\vol(R)\:\sum_{\mskd \in I_n(R)}\:u_{\mskd}\mu_{\mskd}.
\end{equation}
\end{lemma}
Obviously, for each sequence $\{u(n):n\in {\bbN}\}$ weakly converging to some
$u \in W_2^1({\bbR}^d)$ the following equality holds:
$\lim_n\: \lev u(n)| \mu\des = \lev u |\mu\des$.

From this Lemma we have
\[ |\lev{\bf u}_n|\msmu_n\des_R| \ \leq \ h^{-d/2}\,q_R( {\bf u}_n)^{1/2}\,
 \Vert \mu \Vert_{2,-1}.\]

Inequalities (\ref{exp4.4}), (\ref{exp5.5}) and the variational equalities 
(\ref{exp4.1}) imply the first result towards our proof of convergence of 
approximate solutions. If ${\bf u}_n$ solve (\ref{exp4.1}) or (\ref{exp5.3}) then
\begin{equation}
 h^d \,\Md\:\lnorm {\bf u}_n \lnorm_{2,1}^2 \ \leq \ h^d\,
 \lev {\bf u}_n\,|\,(\lambda I+A_n)\,{\bf u}_n \des_R \ \leq \  
 \Vert u(R,n) \Vert _{2,1}\,\Vert \mu \Vert_{2,-1}.
\end{equation}

\begin{corollary}\label{corn5.2} Let ${\bf u}_n = T(\lambda,A_n)\msmu_n$ and
$u(R,n) = \Phi_n(R) {\bf u}_n$. Then
for each $R$ the sequence $\mathfrak{U} = \{ u(R,n): n \in {\bbN}\} 
\subset \cup_n E_n(R,{\bbR}^d)$ converges weakly in $W_2^1({\bbR}^d)$ to some
$u \in W_2^1({\bbR}^d)$.
\end{corollary}

Let $u^\ast$ be the solution of (\ref{exp2.10}). Then the sequence of functions 
$\hat{u}^\ast(n)$, defined by~(\ref{exn3.13}), strongly converges to $u^\ast$ in 
$W_2^1({\bbR}^d)$. In the 
remaining part of this analysis we have to demonstrate the expected property 
$\lim_n u(R,n) = \lim_n \hat{u}^\ast(R,n) = u^\ast$ for each $R$. We follow the 
well-known finite element technique.

\begin{equation}\begin{array}{rcl}
 \Md\, h^d\,\lnorm {\bf u}_n-\hat{{\bf u}}_n^\ast \lnorm_{R2,1}^2 &\leq&
 h^d\,\lev{\bf u}_n-\hat{{\bf u}}_n^\ast \,|\,(\lambda I+A_n)\,
 ({\bf u}_n-\hat{{\bf u}}_n^\ast )\des_R \\  &=&
 h^d\,\lev {\bf u}_n-\hat{{\bf u}}_n^\ast\, |\,(\lambda I+A_n)\, {\bf u}_n\des_R \\  &-&
 h^d\,\lev {\bf u}_n-\hat{{\bf u}}_n^\ast \,|\,(\lambda I+A_n)\,\hat{{\bf u}}_n^\ast \des_R \\
 = h^d\,\ \lev {\bf u}_n-\hat{{\bf u}}_n^\ast\,|\, \msmu_n\des_R  &-&
 h^d\,\lev {\bf u}_n-\hat{{\bf u}}_n^\ast \,|\,(\lambda I+A_n)\,\hat{{\bf u}}_n^\ast \des_R.
 \end{array}
\end{equation}
By Lemma \ref{lem5.2} the first term on the right hand side converges to 
$\lev u-u^\ast |\mu\des$. By the consistency property of Proposition 
\ref{Prop5.1} the second term converges to the same value. 

\begin{theorem}\label{Th5.1} Let $\mathfrak{U}$ be as in Corollary \ref{corn5.2}.
Then the sequence $\mathfrak{U}$ converges $W_2^1({\bbR})$-strongly to the unique 
solution $u^\ast$ to (\ref{exp2.10}).
\end{theorem}

From this result, Lemma \ref{lem5.2} and Lemma \ref{lem3.1} we get another important 
result for $\lambda = 0$.

\begin{corollary}\label{corn5.3} 
Let $D$ be a bounded domain with Lipsithz boundary and $\mu \in W_2^{-1}(D)$.
Let $A_n(D)$ be the restriction to $G_n(R,D)$ of $A_n$, $\msmu_n$ on
$G_n(R,D)$ satisfy (\ref{exp5.5}) and ${\bf u}_n = A_n(D)^{-1}\msmu_n$.
Then the sequence $\mathfrak{U}$ converges 
strongly in $\dot{W}_2^1(D)$ to the unique weak solution $u$ of (\ref{exp2.9}).
\end{corollary}

\section{Convergence in $L_1$-spaces}\label{sec6}
In this section we consider Problem~(\ref{exp2.9}) for a bounded domain 
$D$ and numerical solutions in $\dot{W}_p^1(D)$. 
Furthermore, for the sake of simplicity of exposition, we restrict our
analysis to the differential operator $A_0(\msx) = -\sum_{ij}\partial_i
a_{ij}(\msx)\partial_j$ and its discretizations $A_n(D)$. We consider the 
boundary value problem (\ref{exp2.9}) with $\mu \in {\cal R}(D)$ 
and its discretizations 
\begin{equation}\label{ex6.1}\displaystyle
 A_n(D)\,{\bf u}_n \ = \ \msmu_n
\end{equation}
on $G_n(R,D)$, where $A_n(D)$ have the compartmental structure. The discretizations 
$\msmu_n \in l(G_n(R,D))$ are defined by
\[\mu_{\mskd} \ = \ \frac{\psi_{\mskd}}{\Vert \psi_{\mskd}\Vert_1}\:
 \lev \psi_{\mskd}\,|\, \mu\des.\]
Apparently, $\lev v(n)|\mu\des = h^d\lev {\bf v}_n|\msmu_n\des_R$ for
any $v(n) = \sum v_{\mskd}\psi_{\mskd} \in E_n(R,D)$. Problem (\ref{ex6.1}) is
defined for grid-functions on $G_n(R,D)$. 

\subsection{Boundedness of approximate solutions}
For any pair $r, R, \ 0 < r < R \leq 1$ and $\msv \in D$ we define the balls 
$B_r(\msv), B_R(\msv)$. By using the functions $u \in \dot{W}_2^1(D)$ we define 
the sets $A(r,s,\msv) = \{\msx \in B_r(\msv) \,:\,u(\msx) \geq s\}$, 
where $s \in {\bbR}$. The measure of $A(r,s,\msv)$
is denoted by $a(r,s,\msv) = \mess(A(r,s,\msv))$. 
Let us assume that there exist a subset $\EuScript{G} \subset \dot{W}_2^1(D)$ 
and two numbers, $c_1, c_2$, independent of $r, R, s,\msv$, such that the
inequality
\begin{equation}\label{ex6.3}\displaystyle
 \sum_{i=1}^d\: \Vert \,{\bbJ}_{A(r,s,\msvd)}\,\partial_i u\, \Vert_2^2 \ \leq \
 c_1 \:\Vert\,{\bbJ}_{A(R,s,\msvd)}\,(u-s)\,\Vert_2^2\:\left [
 1\:+\;\frac{1}{(R-r)^2} \right ] \:+\:
 c_2 \:s^2\:\big(a(R,s,\msv)\big)^{1-2/q}
\end{equation}
is valid for a fixed $q > d$, each $u \in \EuScript{G}$, and all $s \geq s_0$ with 
some $s_0$. Then~\cite{LU} the functions $u \in \EuScript{G}$ are bounded on $D$ with a bound 
which depends on $D, c_1,c_2$ and $\Vert u \Vert_2$.
Instead of balls $B_r(\msv)$ one can use rectangles
$S_n(\msp,\msv) = \prod_{i=1}^d [-hp_i+hv_i,hv_i+hp_i]$ as well.
Actually, it is sufficient to consider balls (rectangles) with centers 
$\msv$ in a dense set of $D$ whose radii (edges)
are contained in a sequence $\{r_m: m \in {\bbN}\} \subset (0,1)$,
converging to zero. For instance, the sets $S_n(\msp,\msv)$
and $\msv \in G_n(R),$ match this weaker
condition. This fact enables a straightforward application of
(\ref{ex6.3}) to the sequence of functions in Corollary \ref{corn5.3}.
First, we have to establish a discretized version of~(\ref{ex6.3}), and then
we have to prove that the constructed discretized version implies~(\ref{ex6.3})
for the sequence $\mathfrak{U}$ of Corollary \ref{corn5.3}. 
Further, let
\[\begin{array}{lll}
 A(\msr,s,\msv)  &=& S_m(\msp,\msv) \:\cap\:\{\,u(\msx) \,\geq\,s\,\},\\
 A(\msR,s,\msv)  &=& S_{m-t}(\msp,\msv) \:\cap\:\{\,u(\msx) \,\geq\,s\,\},
 \end{array}\]
where the sets $A(\msR,s,\msv)$ are larger than $A(\msr,s,\msv)$,
${\rm dist}(\partial A(\msr,s,\msv), \partial A(\msR,s,\msv)) = h(m)(2^t-1)
\underline{p}$, where $\underline{p} = \min p_k$ and $m >t$. 
Hence, the symbols $\msr, \msR$ stand for the $d$-dimensional 
parameters $h(m)\msp$ and $h(m-t)\msp$, respectively. The discretization of 
$A(\msr,s,\msv)$ is defined by $F_n(\msr,s,\msv) = G_n(R,A(\msr,s,\msv))$. 
The index set of $F_n(\msr,s,\msv)$ is denoted by $J_n(\msr,s,\msv)$ 
and its cardinal number by $j_n(\msr,s,\msv) = {\rm card}(J_n(\msr,s,\msv))$. 
\begin{figure}
\caption{The function $\chi(\msr,\msR,\msv)$}
\label{fgn6.1}
\beginpicture
\setcoordinatesystem units <0.2pt,0.2pt>
\setplotarea x from -600 to 0, y from -100 to 700
\setlinear
\plot 0 0, 400 0, 400 600, 0 600, 0 0 /
\plot 30 30, 370 30, 370 570, 30 570, 30 30 /
\plot 150 150, 250 150, 250 450, 150 450, 150 150 /
\plot 350 500, 370 510, 350 520 /
\plot 420 500, 400 510, 420 520 /
\plot 30 200, 150 200 /
\plot 50 190, 30 200, 50 210 /
\plot 130 190, 150 200, 130 210 /
\plot 200 200, 450 160 /
\plot 210 187, 200 200, 215 210 /
\plot 300 80, 450 0 /
\plot 315 83, 300 80, 305 65 /
\plot 180 585, 180 670 /
\plot 170 605, 180 585, 190 605 /
\put{$K(\msr)$} [c] at 250 350
\put{$K(\msR)$} [c] at 175 500
\put{$K_+(\msR)$} [c] at 450 300
\put{$p_1 h_n$} [c] at 460 560
\put{$R_1$-$r_1$} [c] at 110 120
\put{"roof"} [c] at 540 160
\put{of $\msch$} [c] at 540 120
\put{"slope"} [c] at 540 20
\put{of $\msch$} [c] at 540 -20
\put{"room" for support} [l] at 220 695
\put{of $Z(i,\pm p_i)\msch$} [l] at 220 650
\endpicture
\end{figure}
\par
Let $r<R \in {\bbN}$ and let us define the cut-off function 
$\theta(r,R,\cdot)$ on ${\bbR}$ 
as a continuous piecewise linear function, such that $\theta(r,R,x) = 1$
for $x \in [-r,r]$, and zero outside of $[-R,R]$. Apparently,
$\theta(r,R,\cdot)' = (1/(R-r)) {\bbJ}_{[-R,-r]} - (1/(R-r)) {\bbJ}_{[r,R]}$. 
Let $\msr = (r_1,r_2,\ldots,r_d)$ and $\msR = (R_1,R_2,\ldots,R_d)$.
The continuous function $\msx \mapsto \chi(\msr,\msR,\msv,\msx)
= \prod_{i=1}^d\theta(r_i,R_i,|x_i-v_i|)$ is a cut-off function on ${\bbR}^d$
for which $|\partial \chi(\msr,\msR,\msv,\cdot)| \leq \max_i
(R_i-r_i)^{-1}$. 

Discretizations of $\chi(\msr,\msR,\msv,\cdot)$ on $G_n(R)$ are denoted
by $\msch_n(\msr,\msR,\msv)$ and they are defined in terms of rougher grids
$G_{m-t} \subset G_m \subset G_n$ for which $m-t < m < n$. For each
$\msv \in G_n(R)$, we define three sets
\begin{equation}\label{exp6.0}\begin{array}{l}
 K(\msr) \ = \ S_{m}(\msp,\msv) = 
 \{\msx : \chi(\msr,\msR,\msv,\msx)=1\},\\
 K(\msR) \ = \ S_{m-t}(\msp,\msv)=
 \supp(\chi(\msr,\msR,\msv)),\\ \displaystyle
 K_+(\msR) \ = \ K(\msR)+S_n(\msp,{\bf 0}) \ = \ \cup_{\msvd \in
 K(\msRd)} S_n(\msp,\msv),\end{array}
\end{equation}
which are illustrated in Figure \ref{fgn6.1}. The rectangles (\ref{exp6.0})
have edges $2r_i= p_i 2^{1-m}, 2R_i=p_i 2^{1+t-m}$ and $p_i (2^{1+t-m}+
2^{1-n})$ units, respectively. It is important to notice that the 
grid-functions $\msch_n(\msr,\msR,\msv)$ are defined by sets $S_n(\msp,\msv)$ 
which are related to the numerical neighbourhoods of constructed schemes.
In our next proof we again use the fact that the sets $\supp(\max\{
a_{ij},0\}), i \ne j$ are strictly separated, so that we can assume that the
functions $a_{ij}$ on rectangles (\ref{exp6.0}) do not change sign.

It is easy to verify that
$|U_i(p_i) \msch_n(\msr,\msR,\msv)|\leq (h(m)\underline{p}(2^t-1))^{-1}$, 
so that we can write
\[|U_i(p_i)\,\msch_n(\msr,\msR,\msv)| \ \leq \ \frac{1}{\underline{p}^2}\:
 \frac{1}{h(n)^2\,(L-l)^2} \ \leq \ \frac{\rho}{|\msR-\msr|},\]
where $L=2^{t+n-m}, l =2^{n-m}$, and $\rho$ is a number depending on $\underline{p}
=\min p_i(l)$.
The sets $F_n(\msr,s,\msv)$ and $F_n(\msR,s,\msv)$ are defined in terms of
the sets (\ref{exp6.0}),
\[\begin{array}{l}
 F_n(\msr,s,\msv) \ = \ K(\msr)\cap \supp({\bf w}_n),\\
 F_n(\msR,s,\msv) \ = \ K_+(\msR)\cap \supp({\bf w}_n).\end{array}\]
The following estimates are used in our next proof:
\begin{equation}\label{exp6.10}\begin{array}{l}\displaystyle
 \max_{l,i} \ Z(i,\pm p_i(l))\chi(\msr,\msR,\msv) \ \leq \ {\bbJ}_{F(\msRd,s,\msvd)},
 \\ \displaystyle
 \max_{l,i} \ \big|U_i(p_i(l))\,\big(Z(i,\pm p_i(l))\msch_n(\msr,\msR,\msv)\big)\big|
 \ \leq \ \frac{\rho}{|\msR-\msr|} {\bbJ}_{F(\msRd,s,\msvd)}.\end{array}
\end{equation}
Now we consider the following auxiliary problem:
\begin{equation}\label{ex6.4}
 A_n(D)\,{\bf u}_n \ = \ {\bf f}_{0n}\,-\,\sum_{i=1}^d\:V_i(p_i)\,{\bf f}_{in}
 \ := \ \msmu_n,
\end{equation}
where ${\bf u}_n$ are grid-solutions and ${\bf f}_i$ are grid-functions on 
$G_n(R,D)$.
\begin{lemma}\label{lem6.1} Let ${\bf u}_n = A_n(D)^{-1} \msmu_n$, where $\msmu_n$ are
defined by (\ref{ex6.4}). There exist positive numbers $c_1, c_2$ and $s_0$, 
independent of $n, \msv$ such that the inequalities
\begin{equation}\label{ex6.5}\begin{array}{lll}\displaystyle
 \sum_{i=1}^d\:\lnorm \msch_n(\msr,\msR,\msv)\,U_i(p_i){\bf u} \lnorm_{R2}^2
 &\leq&  \displaystyle  c_1\:\lnorm {\bbJ}_{F(\msRd,s,\msvd)} ({\bf u}_n -s) \lnorm_{R2}^2
 \:\left [ 1 \,+\, \frac{1}{|\msR-\msr|^2}\right ] \\
 &+&\displaystyle
 c_2 \:s^2\:\sum_{i=0}^d\:\lnorm {\bbJ}_{F(\msRd,s,\msvd)} {\bf f}_{in}\lnorm_{Rq}^2
 \big(j(\msR,s,\msv)\big)^{1-2/q}\end{array}
\end{equation}
are valid for $s \geq s_0$ and $q > 2$.
\end{lemma}
{\Proof} Because of the rule~(\ref{exp4.14}) of construction of discretizations 
for $d > 2$, it is sufficient to prove~(\ref{ex6.5}) for two-dimensional
problems. Basic schemes are considered first. 

The function $w = \max\{u -s,0\}$ has the discretizations ${\bf w}_n$ 
and ${\bf w}_nU_i {\bf w}_n = {\bf w}_nU_i {\bf u}_n$. 
Now we have to evaluate the forms
\begin{equation}\label{exp6.2}\begin{array}{l}
 a_n^{(-)}(\chi^2 w,u) \:=\:
 \lev \msch_n(\msr,\msR,\msv)^2
 {\bf w}_n\,|\, A_n^{(-)} \,{\bf u}_n\des_{R},\\
 a_n^{(+)}(\chi^2 w,u) \:=\:
 \lev \msch_n(\msr,\msR,\msv)^2
 {\bf w}_n\,|\, A_n^{(+)} \,{\bf u}_n\des_{R},\end{array}
\end{equation}
in terms of matrices $\EuScript{A}_n^{(\pm)}(i,j,r_i,r_j)$ of Expression 
(\ref{exp4.19}). It is sufficient to consider the form $a_n^{(-)}$ and one
of the terms. For instance, the term involving the matrix $A_n^{(-)}(l)$ is 
$-\sum_{ij =1}^2\:V_i(p_i)\EuScript{A}_n^{(+)}(i,j,p_i,p_j)\,U_j(p_j)$. 

In order to write expressions in this proof in a concise form we omit
various indices in the notation. Thus we use ${\bf u}, {\bf w},
\msch, U_i$ and $\EuScript{A}^{(\pm)}(i,j)$. Apart from this, we need the
notations $Z_i{\bf f}=Z_n(p_i,i){\bf f}_n$, so that we can write
$U_i{\bf f}_n = (hp_i)^{-1}(Z_i{\bf f}-{\bf f})$. In the equality
$U_i({\bf g}{\bf w}) = (U_i {\bf g}) Z_i{\bf w} + {\bf g}U_i {\bf w}$ we 
insert ${\bf g} = \msch^2$ to get:
\begin{equation}\label{exp6.6}
 U_i\big(\msch^2 {\bf w}\big) \ = \
 \msch^2\,U_i{\bf w} \,+\, \big(U_i\msch\big)\,\big(Z_i\msch +\msch
 \big) \,Z_i{\bf w}.
\end{equation}
Also we need
\begin{equation}\label{exp6.3}
 \big(Z_i\msch +\msch \big) \,U_j{\bf w} \ = \ 2 \msch U_j{\bf w}+
 \frac{p_i}{p_j}(U_k \msch) Z_j{\bf w} - \frac{p_i}{p_j}(U_k \msch){\bf w}.
\end{equation}
In this way the finite differences $U_j{\bf w}$ are always multiplied by
the grid-function $\msch$ (not with $Z_k\msch$).

We start the proof by an analysis of~(\ref{exp6.2}):
\[\begin{array}{l} \displaystyle
 \lev \msch^2\, {\bf w} \,|\, A_n^{(-)}(l) \,{\bf u}\des_{R} \ = \
 \sum_{i,j=1} \lev \ \msch\,  U_i{\bf w}\,|\, \EuScript{A}^{(+)}(i,j) \,\msch\,
 U_j \,{\bf w}\des_{R} \\ \displaystyle
 +\,  \sum_{i,j=1} \lev (Z_i{\bf w})\, U_i\msch\,|\, 
 \EuScript{A}^{(+)}(i,j) \,\big(Z_i\msch+\msch\big)\,U_j \,
 {\bf w}\des_{R}. \end{array}\]
The first term on the right hand side can be estimated from bellow as a consequence of the
strict ellipticity~(\ref{exp4.4}) of Proposition~\ref{Prop4.2}. 
The left hand side of this equality is equal to 
$\lev \msch^2 {\bf w},| \msmu \des_{R}$, so that we can write:
\begin{equation}\label{exp6.7}\begin{array}{l} \displaystyle
 \Md \:\sum_{i=1}\: \lnorm \msch \,U_i {\bf w}  \lnorm_{R2}^2 \ \leq \ 
 \lev \msch^2\, {\bf w} \,|\, \msmu\des_{R} \\ \displaystyle
 -\: \sum_{i,j=1} \lev Z_i{\bf w}\, U_i\msch \,|\, 
 \EuScript{A}^{(+)}(i,j) \,\big(Z_i\msch+\msch\big)\,U_j 
 \,{\bf w}\des_{R}. \end{array}
\end{equation}
The two terms on the right hand side are denoted by $T_i, i = 1,2$, respectively.
Let us estimate them from above. The representation of $\msmu$, in~(\ref{ex6.4}), 
is inserted into the first term to get:
\[ |T_1| \ = \ |\lev \msch^2\,{\bf w} \,|\, \msmu\des_{R}| \ \leq \ 
 |\lev \msch^2\, {\bf w}|{\bf f}_0\des_{R}| \,+\,\sum_{i=1}^2\:
 |\lev U_i(\msch^2\,{\bf w})|{\bf f}_i\des_{R}|.\]
First (\ref{exp6.6}) and then (\ref{exp6.3}) are utilized in the last two 
terms on right hand side.
The functions $\msch, Z_k\msch$ are estimated from above by~(\ref{exp6.10})
and then the CSB inequality and the inequality $a b \leq (\varepsilon/2) a^ 2 + 
(1/2\varepsilon) b^2$ with a convenient choice of $\varepsilon$, is applied to all 
the terms on the right hand side. 
\[\begin{array}{lll} \displaystyle
 |T_1| &\leq&  \displaystyle
 \varepsilon\,\sum_{i=1}^2\:\lnorm \msch \,U_i{\bf w} \lnorm_{R2}^2
 \,+\,b_1\Big (1\,+\,\frac{1}{|\msR-\msr|^2}\Big)\:
 \lnorm {\bbJ}_{F(\msRd,s,\msvd)}\,{\bf w}\lnorm_{R2}^2 \\
 &+& \displaystyle
 \left(1+\frac{1}{4\varepsilon}\right)\,\sum_{i=0}^2\:
 \lnorm {\bbJ}_{F(\msRd,s,\msvd)}\,{\bf f}_i\lnorm_{R2}^2 
 ,\end{array}\]
where $b_1$ does not depend on $n, \varepsilon$. It is important to
notice that the first term on the right hand side of this inequality and the left
hand side of~(\ref{exp6.7}) involve the same norms, and can therefore be
subtracted. Now we use the H\"{o}lder
inequality $\lnorm {\bbJ}_{F(\msRd,s,\msvd)} {\bf f}_i\lnorm_{R2}^2 \leq \lnorm 
{\bbJ}_{F(\msRd,s,\msvd)} {\bf f}_i\lnorm_{Rq}^2 
j(\msR,s,\msv)^{1-2/q}$ and get the first intermediate result:
\[\begin{array}{lll} \displaystyle
 |T_1| &\leq&  \displaystyle
 \varepsilon\,\sum_{i=1}^2\:\lnorm \msch \,U_i{\bf w} \lnorm_{R2}^2
 \,+\,b_1\:\Big (1\,+\,\frac{1}{|\msR-\msr|^2}\Big)\:
 \lnorm {\bbJ}_{F(\msRd,s,\msvd)}\,{\bf w}\lnorm_{R2}^2 \\
 &+& \displaystyle
 b_2\,\sum_{i=0}^2\: \lnorm {\bbJ}_{F(\msRd,s,\msvd)}\,{\bf f}_i\lnorm_{Rq}^2 
 \:j(\msR,s,\msv)^{1-2/q},\end{array}\]
where $b_2$ depends on $1/\varepsilon$.

The quantity $T_2$ can be estimated analogously. Let us use (\ref{exp6.3})
in Expression (\ref{exp6.7}) and apply the same technique as in the previous 
case. We get:
\[ |T_2| \ \leq \ \varepsilon\,\sum_{i=1}^2\:\lnorm \msch \,U_i{\bf w} \lnorm_{R2}^2
\,+\,b_3\:\frac{1}{|\msR-\msr|^2}\:
 \lnorm {\bbJ}_{F(\msRd,s,\msvd)}\,{\bf w}\lnorm_{R2}^2,\]
where $b_3$ depends on $\Mg,1/\varepsilon$.
Thus the inequality (\ref{exp6.7}) implies the following result:
\[\begin{array}{lll}\displaystyle
 \big(\Md \!-\!  2\varepsilon ) \:\sum_{i=1}^2\: \lnorm \msch\,U_i {\bf u}
 \lnorm_{R2}^2  
 &\leq& \displaystyle
 b_4\:\Big (1\,+\,\frac{1}{|\msR-\msr|^2}\Big)\:
 \lnorm {\bbJ}_{F(\msRd,s,\msvd)}\,{\bf w}\lnorm_{R2}^2 \\
 &+& \displaystyle
 b_2\,\sum_{i=0}^2\: \lnorm {\bbJ}_{F(\msRd,s,\msvd)}\,{\bf f}_i\lnorm_{Rq}^2 
 \:j(\msR,s,\msv)^{1-2/q},\end{array}\]
which is equivalent to the assertion of lemma. 

Now we consider the extended schemes. Let the discretizations $A_n(D)$ of
(\ref{ex6.5}) be constructed by extended schemes and let ${\bf u}_n$ be the 
corresponding solution. Then we have $a_n(v,u)=\lev {\bf v}_n|\msmu_n\des$,
where $a_n(v,u)$ is defined in~(\ref{exp4.9}). Let us define ${\bf v}_n
=\msch^2 {\bf w}$ and let its image $v(n)=\Phi_n {\bf v}_n$ be
inserted into the expression $a_n(v,u)$. If we apply (\ref{exp6.6}) to
the function $v(n)$ we get:
\begin{equation}\label{exp6.11}
 a_n(v(n),u) \ = \ a_n(\chi;w,w)+\sum_s b_n\big(s,\msch,{\bf w}\big),
\end{equation}
where the terms on the right hand side are defined as follows. The form
$a_n(\chi;w,w)$ is the sum of forms (\ref{exp5.16}) and (\ref{exp5.17}) in 
which the factors $a_{ij}(\msx+\mst^{(\pm)}(\msr))$ are replaced with the 
factors $\msch(\msx)^2 a_{ij}(\msx+\mst^{(\pm)}(\msr))$. Each
$b_n(s,\msch,{\bf w})$ is a form which is a first degree polynomial in
$U_i(p_i(l))){\bf w}_n$, with various $l,i,p_i(l)$, resulting from
the application of the rule~(\ref{exp6.6}). 
As in the corresponding proof of Proposition~\ref{Prop4.2}, the quantity
$a_n(\chi;w,w)$ can be estimated from bellow,
\[ a_n(\chi;w,w) \ \geq \ \omega(a)\,\sum_{i=1}^2\: \lnorm \msch\,U_i {\bf w}
 \lnorm_{2}^2 .\]
An attempt to estimate the remaining terms on the right hand side in~ (\ref{exp6.11})
from above as in the first part of proof would lead to an
unresolvable problem. Among the resulting terms there would be
$\varepsilon \sum_i \lnorm \msch_n U_i{\bf w}_n\lnorm_2^2$ as well as
$\varepsilon \sum_i \lnorm \msch_n U_i(p_i){\bf w}_n\lnorm_2^2$. The former
one can be moved to the left hand side of (\ref{exp6.11}) and subtracted from
$\omega(a)\sum_{i} \lnorm \msch_n U_i {\bf w}_n\lnorm_{2}^2$. Unfortunately,
this cannot be done with the latter one. Therefore, before finding upper bounds 
on $b_n(s,\msch,{\bf w})$ we have to apply another version of~(\ref{exp6.3}):
\[\begin{array}{l}
 \big(\msch_n +Z_n(p_i,i)\msch_n\big)U_j(p_j) {\bf w}_n \ = \ 
 2 \tilde{\msch}_n U_j(p_j) {\bf w}_n\,+ \\ \displaystyle
 \frac{p_i}{p_j}\big(U_i(p_i)\tilde{\msch}_n\big) \big(Z_n(p_j,j)-I\big)
 {\bf w}_n +
 \frac{1}{p_j}\Big[\big(I+Z_n(p_i,i)\big)\big(\msch_n-\tilde{\msch}_n\big)
 \Big] \big(Z_n(p_j,j)-I\big) {\bf w}_n. \end{array}\]
If these expressions are applied, then estimates from above of $b_n(s,\msch,{\bf w})$
contain the term $\varepsilon \sum_i \lnorm \tilde{\msch}_n U_i(p_i)
{\bf w}_n\lnorm_2^2$. Let us choose $\tilde{\chi}(\msx) = \chi(\msr,\msR,
\msv,\msx)\chi(\msr,\msR,\msv,\msx+h(p_1\mse_1+_2\mse_2))$. Then
\[ \supp(\tilde{\chi}) +C_n(\msp,{\bf 0}) \ = \ K(\msR),\]
and $(\tilde{\msch}_n U_j(p_j) {\bf w}_n)_{\mskd}$ is a linear combination of
$(U_j(1) {\bf w}_n)_{\msld}$ such that $h\msl \in K(R)$. Now we have
\[ \varepsilon \sum_i \lnorm \tilde{\msch}_n U_i(p_i){\bf w}_n\lnorm_2^2
 \ \leq \ \varepsilon \sum_i \lnorm \msch_n U_i(1){\bf w}_n\lnorm_2^2 \]
and the difficulty regarding $\lnorm \msch_n U_i(p_i){\bf w}_n\lnorm_2^2$
with various $p_i \ne 1$ is removed. The grid-functions $h_{-1}(\msch_n-\tilde{\msch}_n)$ 
have estimates from above as in~(\ref{exp6.10}). In this way we get (\ref{ex6.5}) again. {\QED}

The function $w(n)=\max\{u(n)-s,0\}$ and the grid-function ${\bf w}_n$ which is 
componentwisely by $w_{\mskd}=\max\{u_{\mskd}-s,0\}$, are related by the
expression $w(n) = \Phi_n(R) {\bf w}_n$. The index set of grid-knots of $G_n(R)$ in
the (closed) set $S_m(\msp,\msv)$ is denoted by $I(m,\msp,\msv)$. The restriction of $w(n)$ 
to the set $S_m(\msp,\msv)$ and the function $w(n,I) = \sum_{\mskd \in I(m,\mspd,\msvd)}
w_{\mskd}\psi_{\mskd}$ are not equal. They are equal on the set $S_m(\msp,\msv)$
and the following rough estimate is valid elsewhere:
\[\begin{array}{lll}\displaystyle
 \Vert w(n,I) \Vert_2^2 &=& \displaystyle
 \Vert {\bbJ}_{S_m(\mspd,\msvd)}\,w(n,I) \Vert_2^2
 \:+\:\Vert (1-{\bbJ}_{S_m(\mspd,\msvd)})\, w(n,I) \Vert_2^2  \\ &\leq& \displaystyle
 2^d\,\Vert {\bbJ}_{S_m(\mspd,\msvd)}\,w(n,I) \Vert_2^2 \:=\:
 2^d\,\Vert {\bbJ}_{S_m(\mspd,\msvd)}\,w(n) \Vert_2^2.\end{array}\]
Now we combine the result of Theorem \ref{thn3.1}, {\em i.e.}:
\[ (1-\sigma^2)\,\lnorm {\bbJ}_{I(m,\mspd,\msvd)}\,{\bf w}_n \lnorm_{R2}^2
 \ \leq \ \Vert w(n,I)\Vert_2^2 \ \leq \ \lnorm {\bbJ}_{I(m,\mspd,\msvd)}\,
 {\bf w}_n \lnorm_{R2}^2 \]
and the previous inequality in order to get:
\begin{equation}\label{exp6.1}
 \lnorm {\bbJ}_{I(m,\mspd,\msvd)}\,{\bf w}_n \lnorm_{R2}^2 \ \leq \
 \frac{1}{1-\sigma^2}\,\Vert w(n,I) \Vert_2^2 \ \leq \
 \frac{2^d}{1-\sigma^2}\,\Vert {\bbJ}_{S_m(\mspd,\msvd)}\,w(n) \Vert_2^2
\end{equation}
It remains to compare the functions $\partial_iw(n)$ and $U_i(p_i){\bf w}_n$.
By Expression~(\ref{exp6.8}) we easily prove:
\begin{equation}\label{exp6.13}\begin{array}{c}
 \Vert {\bbJ}_{S_m(\mspd,\msvd)}\,\partial_iw(n) \Vert_2^2 \ \leq \ 
 \Vert \sum_{\mskd \in I(m,\mspd,\msvd)}\:\partial_i w(n)
 \Vert_2^2 \ \leq \ h^d\,\lnorm {\bbJ}_{I(m,\mspd,\msvd)}\,
 \big(U_i(p_i){\bf w}_n\big)_{\mskd} \lnorm_{R2}^2\\
 \ \leq \ h^d\,\lnorm \msch_n(\msr,\msR,\msv)\,U_i(p_i){\bf w}_n\lnorm_{R2}^2.
 \end{array}
\end{equation}
Let us multiply (\ref{ex6.5}) by $h^d$, use (\ref{exp6.13}) and afterwards 
(\ref{exp6.1}). The result is:
\begin{equation}\label{exp6.14}\begin{array}{lll}
 \sum_{i=1}^d\:\Vert {\bbJ}_{A(\msrd,s,\msvd)}\,\partial_i w(n) \Vert_2^2 &\leq&
 \displaystyle
 \frac{2^d\,c_1}{1-\sigma^2}\,\Vert {\bbJ}_{A(\msRd,s,\msvd)}\,w(n) \Vert_2^2
 \,\Big(1+\frac{1}{|\msR-\msr|}\Big) \\
 &+& \displaystyle
 c_2\,s^2\,\big(h^{2d/q}\,\sum_{i=0}^d\,\lnorm {\bf f}_{in}\lnorm_{Rq}^2\big)\:
 a(\msR,s,\msv)^{1-2/q}.\end{array}
\end{equation}

In this way we obtain an important property of grid-solutions to Problems~(\ref{ex6.4}).

\begin{corollary}\label{cor6.1} Let $D \subset {\bbR}^d$ be a bounded domain with 
a Lipshitz boundary and $q > d$. If $h^d \sum_{i=0}^d\lnorm {\bf f}_{in}\lnorm_{Rq}^q$
are bounded by a number uniformly with respect to $n$, then the grid-solutions
${\bf u}_n$ on $G_n(R,D)$ to (\ref{ex6.4}) are also bounded uniformly with 
respect to $n$.
\end{corollary}
{\Proof} One has to compare (\ref{ex6.3}) and inequality~(\ref{exp6.14}). {\QED}
\subsection{Weak consistency and convergence}
As in the case of $W_2^1$-approach, the convergence in $L_p$-spaces is
proved by utilizing certain kind of consistency. This weaker consistency
is defined in terms of sequences similar to~(\ref{exp5.7}):
\begin{equation}\label{exp6.9}\begin{array}{lllll} \displaystyle
 \mathfrak{V}(p,c) &=&  \{v(n)\::\: n \in {\bbN}, \ h^d\,\sum_{i}\:
 \lnorm U_i{\bf v}_n\lnorm_p^p 
 \:\leq \:c \} &\subset& \cup_n E_n({\bbR}^d),\\
 \mathfrak{U}_0 &=&  \{\hat{u}(n)\::\: n \in {\bbN}, \ u \in C_0^{(1)}(
 {\bbR}^d) \} &\subset& \cup_n E_n({\bbR}^d).
 \end{array}
\end{equation}

\begin{definition}[Weak consistency]\label{def6.1} We say that forms $a_n(\cdot,\cdot)$
on $E_n({\bbR}^d) \times E_n({\bbR}^d)$ are weakly consistent with the form 
(\ref{exp2.4}) if
\[ a(v,u) \ = \ \lim_n h^d\,a_n(v(n),\hat{u}(n)) \]
is valid for any weakly convergent $\mathfrak{V}(p,c), p \in [1,\infty), c > 0$ and any $\mathfrak{U}_0$ 
of (\ref{exp6.9}).
\end{definition}

\begin{lemma}\label{lem6.3} Let a sequence of matrices $\{A_n:n \in {\bbN}\}$
be constructed by basic or extended schemes, and let $\{a_n(\cdot,\cdot):
n \in {\bbN}\}$ be the corresponding sequence of discretizations of (\ref{exp2.4}).
Then the forms $a_n$ are weakly consistent with the form (\ref{exp2.4}).
\end{lemma}

{\Proof} The present proof and proof of Lemma \ref{lem5.4} are the same up to 
the equality (\ref{exp5.11}). The obtained equality must be estimated now 
by the H\"{o}lder inequality:
\[\begin{array}{c} \displaystyle
 \Big|\gamma_n(l,v(n),u(n))\,-\,\theta_n(l,v(n),u(n))\Big| \ \leq \\ \displaystyle
 \overline{p}\:h^{d/p}\lnorm U_i(r_i) {\bf v}_n \lnorm_{Rp} \:h^{d/q}\,
 \max \{\lnorm \big(Z_n(r_j,j)-I\big) U_j(r_j)\hat{{\bf u}}_{n}\lnorm_{Rq} \,:\,
 j = 1,2,\ldots,d\}. \end{array}\]
Thus we have to consider $\lnorm \big(Z_n(r_j,j)-I\big) U_j(r_j)\hat{{\bf u}}_{n}
\lnorm_{Rq}
$ for a fixed $j \in {\bbN}$ and prove that this quantity converges to zero
as $n \to \infty$. It is easy to get the following expression:
\[ \Big(\big(Z_n(r_j,j)-I\big)\, U_j(r_j)\hat{{\bf u}}_{n}\Big)_{\mskd} \ = \ 
 \frac{1}{\Vert \psi_{\mskd} \Vert_1}\, 
 \big( \psi_{\mskd}\,|\,(Z(th\mse_j)-I\big) \der_j(r_jh)u
 \big) \ = \ \big(\hat{{\bf w}}_n(r_j)\big)_{\mskd},\]
where the function $w(r_j) = (Z(th\mse_j)-I\big) \der_j(r_jh)u$ is continuous and
converges pointwisely to zero as $h \to 0$. By (ii) of Lemma \ref{lem3.7} we get
\[ h^{d/q}\:\lnorm \hat{{\bf w}}_n(r_j) \lnorm_{Rq} 
 \ \leq \ \Vert w(r_j) \Vert_q \ 
 \small{\begin{array}{c} \\ \longrightarrow \\ h \to 0\end{array}}
   \ 0.\]
{\QED}
\begin{theorem}\label{thn6.1} Let $D \subset {\bbR}^d$ be a bounded domain with 
a Lipshitz boundary, and let $\mathfrak{U} = \{ u(n): n \in {\bbN}\}\subset E_n(R,D)$
be a sequence of approximate solutions, $u(n) = \Phi_n(R){\bf u}_n$, where ${\bf u}_n$
are grid-solutions to (\ref{ex6.1}). For each $p \in [1,d/(d-1))$ there exists 
a subset $J(p)\subset {\bbN}$ such that $\mathfrak{U} = \{ u(n): n \in J(p)\}$ 
converges to the unique solution to (\ref{exp2.9}) strongly in $L_p(D)$ and weakly 
in $\dot{W}_p^1(D)$.
\end{theorem}
{\Proof} The proof is split into two steps. In the first step we have to prove 
that for each $p \in (1,d/(d-1))$ there exits $u \in \dot{W}_p^1(D)$ and a subset 
$J(p)\subset {\bbN}$ such that
\[\begin{array}{lllll}
 u &=& \displaystyle w \!\! - \!\! \lim_{n \in J(p)} u(n) &\in& \dot{W}_p^1(D),\\
 u &=& \displaystyle s \!\! - \!\! \lim_{n \in J(p)} u(n) &\in& L_p(D).\end{array}\]
In the second step of the proof we demonstrate that $u$ coincides with the solution
$u^\ast$ of the original problem~(\ref{exp2.9}).

Because of $D \subset S_1(\mst,\msv)$, for some $\mst,\msv$ we have
$\Vert u \Vert_p \leq (\max_k t_k)\,(\max_k \Vert \partial_k u\Vert_p)$ so that
we can consider only $\partial_1 u(n)$ and the corresponding $U_1(p_1){\bf u}_n$.
By~(\ref{exp6.8}) and estimating the function $\psi_{\mskd'}$ by the indicator of its support we get:
\[ |\partial_1 u(n)| \ \leq \ 2^{d}\:\sum_{\mskd \in I_n(R)}\: \Big|\big(U_1(p_1)
 {\bf u}_n)_{\mskd}\Big|\, {\bbJ}_{C_n(\msp,\msv)},\]
where $\msv = h\msk$. Hence
\[ \Vert \partial_1 u(n)\Vert_p^p \ \leq \ 2^{dp}\: h^d\:\lnorm U_1(p_1)
 {\bf u}_n \lnorm_{Rp}^p.\]
It remains to demonstrate that the right hand side is bounded uniformly with
respect to $n$. We can write
\[ \lnorm U_1(p_1){\bf u}_n \lnorm_{Rp} \ = \ \lev {\bf v}_n \,|\,U_1(p_1)\,
 {\bf u}_n \des_{R},\]
where $v_{\mskd} = h^{-d/q}{\rm sign}(U_1(p_1)({\bf u}_n)_{\mskd})|U_1(p_1)
({\bf u}_n)_{\mskd}|^{p-1}/\lnorm U_1(p_1){\bf u}_n \lnorm_{Rp}^{p-1}$. Then
$h^{d/q}\lnorm {\bf v}_n\lnorm_{Rq}=1$. Therefore,
\[ \lnorm U_1(p_1){\bf u}_n \lnorm_{Rp} \ = \ -\,\lev V_1(p_1){\bf v}_n\,|\,{\bf u}_n\des_{R}
 \ = \ -\,\lev A_n^{-1}\, \big(V_1(p_1){\bf v}_n\big)\,|\,\msmu_n\des_{R}. \]
By Corollary \ref{cor6.1} we get $\lnorm A_n^{-1} V_1(p_1){\bf v}_n
\lnorm_{R\infty} \leq \beta$ with some $\beta$ depending on $D$ and the bounds
of $h^d\lnorm {\bf v}_n\lnorm_{Rq}^q$. The final result is $h^{d/p}\lnorm U_i(p_i)
{\bf u}_n\lnorm_{Rp}\leq \beta$. 

The unique solution to (\ref{exp2.9}) is denoted by $u^\ast$. For $p\in (1,d/(d-1))$,
$v \in\dot{C}^{(1)}(\overline{D})$ and the corresponding $\hat{{\bf v}}_n, \hat{v}(n)
=\Phi_n(R)\hat{{\bf v}}_n$, the following relations are valid:
\[\begin{array}{lll}
 1)&\quad& \displaystyle
 v \ = \ s \!\! - \!\! \lim_n \hat{v}(n) \ \in \ \dot{W}_q^1(D),\qquad
 1/p+1/q=1,\\ \displaystyle
 2)&\quad& \displaystyle
 h^d\,\lev \hat{{\bf v}}_n\,|\, \msmu_n\des_{R} \:=\:\lev \hat{v}(n)\,|\,\mu\des 
 \ \rightarrow \ \lev v\,|\, \mu\des,\\ \displaystyle
 3)&\quad& \displaystyle
 \lim_n a(\hat{v}(n),u(n)) \,-\,\lim_n \,h^d\,\lev \hat{{\bf v}}_n\,|\,
 A_n(D){\bf u}_n\des_{R} \ = \ 0. \end{array}\]
The inclusion 1) follows from Lemma~\ref{lem3.6} and the uniform continuity of 
$\partial v_1$ on $\overline{D}$. Although this inclusion is plausible, some technical
details are needed in order to transfer $\partial_1$ from $\partial_1 \hat{v}(n)$
to $\widehat{(\partial_1v)}(n)$. The relation 2) follows from $v \in 
\dot{C}(\overline{D})$ and the construction of $\msmu_n$.
The identity 3) is the weak consistency which is proved in Lemma \ref{lem6.3}. 
Notice that the pair $\hat{v},u$ in this proof and the pair $v,\hat{u}$ of
Lemma~\ref{lem6.3} have  roles interchanged. Now we have
\[\begin{array}{lll}
 a(v,u) &=& \displaystyle a(v,u)\,-\,a(\hat{v}(n),u(n))\\
        &+& \displaystyle a(\hat{v}(n),u(n)) \,-\,h^d\,\lev \hat{{\bf v}}_n\,|\,A_n(D)\,
		{\bf u}_n\des_{R}\\
        &+& \displaystyle h^d\,\lev \hat{{\bf v}}_n\,|\,A_n(D)\,{\bf u}_n\des_{R} \,-\,
 h^d\,\lev \hat{{\bf v}}_n\,|\, \msmu_n \des_{R} \\
        &+& \displaystyle h^d\,\lev \hat{{\bf v}}_n\,|\, \msmu_n \des_{R}. \end{array}\]
The right hand sides in the first three rows are either zero or converge to zero
as $n \to \infty$. Thus $a(v,u) = \lev v|\mu\des$. Because $\dot{C}^{(1)}(\overline{D})$
is dense in $\dot{W}_p^1(D), p \in (1,d/(d-1))$ we get $u = u^\ast$. However $u-u^\ast
=0\in \dot{W}_p^1(D)$ implies $u-u^\ast =0\in \dot{W}_1^1(D)$. {\QED}

\section{Numerical examples in two dimensional case}\label{sec7}
We test the numerical methods on
examples for which we know solutions in the closed form. In the first example
we solve the problem~(\ref{exp2.3}) with $a_{12}=0$ and a 
measure on the right hand side, and in the second exmple we have a piecewise constant
$a_{12}$ and a non-trivial application of the extended scheme.
\begin{figure}
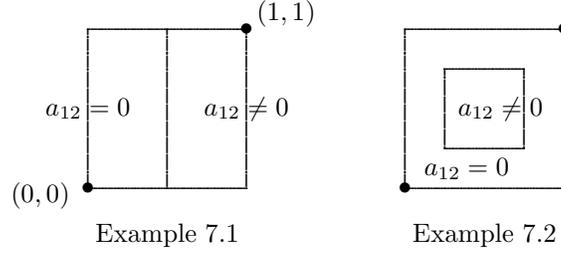

\caption{Domains for Examples~\ref{exam7.1} and~\ref{exam7.2}}
\label{fgn7.1}
\beginpicture
\setcoordinatesystem units <0.6pt,0.6pt>
\setplotarea x from -100 to 400, y from -30 to 110
\setlinear
\plot 50 0, 150 0, 150 100, 50 100, 50 0 /
\plot 100 0, 100 100 /
\plot 250 0, 350 0, 350 100, 250 100, 250 0 /
\plot 275 25, 325 25, 325 75, 275 75, 275 25 /
\put{$\bullet$} [c] at 50 0 
\put{$\bullet$} [c] at 150 100 
\put{$\bullet$} [c] at 250 0 
\put{$\bullet$} [c] at 350 100 
\put{$(0,0)$} [c] at 20 -5
\put{$(1,1)$} [c] at 175 110
\put{$a_{12} \ne 0$} [c] at 150 50
\put{$a_{12} = 0$} [c] at 50 50
\put{Example \ref{exam7.1}} [c] at 100 -30
\put{$a_{12} \ne 0$} [c] at 310 50
\put{$a_{12} = 0$} [l] at 262 12
\put{Example \ref{exam7.2}} [c] at 300 -30
\endpicture
\end{figure}
\begin{example}\label{exam7.1}{\rm
Let the matrix-valued function $\msx \mapsto a(\msx)$ be defined on ${\bbR}^2$ by
\[ a(\msx) \ = \ 10^{-4} \times \left \{\begin{array}{lll}
\left[\begin{array}{cc} 1&0 \\ 0&1\end{array} \right ] & {\rm for} & x_1 < 0.5,
 \\ \left[\begin{array}{cc} \sigma^2&0 \\ 0 &1\end{array}
 \right ] & {\rm for} & x_1 \geq 0.5,\end{array}\right . \]
and let us consider the elliptic differential operator $A(\msx) = -\sum_{ij = 1}^2
\partial_i a_{ij}(\msx)\partial_j$ on ${\bbR}^2$. For $\mst = (0.5,0.5)$ the
differential equation $A(\msx)F(\msx) = \delta(\msx-\mst)$ on ${\bbR}^2$ 
has a solution 
\[ F(\msx) \ = \ -\,c\,\left\{\begin{array}{lll} 
 \ln \big( (x_1-0.5)^2 +(x_2-0.5)^2\big) & {\rm for} & x_1 \:<\:0.5,\\
 \Big[\ln \big( (x_1-0.5)^2 )+\sigma^2(x_2-0.5)^2\big)-\ln 
 (\sigma^2)\Big]
 & {\rm for} & x_1 \:\geq\:0.5,\end{array}\right .\]
where
\[ c \ = \ \Big(2 \pi \,\big(1+{\sigma}\big)\Big)^{-1}.\]
The function $F(\cdot)$ is continuous on ${\bbR}^2 \setminus \{\mst\}$ and its 
first partial derivatives have a jump at $x_1 = 1/2$.
Let $D = (0,1)^2 \subset {\bbR}^2$ as illustrated in Figure~\ref{fgn7.1}.
In this example we consider the boundary 
value problem $A(\msx)u(\msx) = \delta(\msx-\mst)$ with the nonhomogeneous boundary
conditions $u|\partial D = F|\partial D$. Obviously, the solution $u$
coincides with $F$ on $D$.
}\end{example}

In numerical calculations we used $\sigma^2 = 10$. Hence, we have
two sets for the discretization procedure~\ref{dsp4.1}, 
$D_1 = \{\msx \in {\bbR}^2 : x_1 < 0\}$ and
$D_2 = \{\msx \in {\bbR}^2 : x_2 \geq 0\}$. The extended scheme was
utilized for the discretization of $A(\msx)$. 
Discretization of ${\bbR}^2$ is
realized by grid-knots $\msx = hk_1\mse_1 + hk_2\mse_2$, $h=1/400$,
so that $D$ is discretized by $399 \times 399$ grid-knots. The system 
matrix $A_n$ of (\ref{exp5.3}) has the order $399$ and the
linear system is solved iteratively: if $K_n = {\rm diag}(A_n)$
and $Q_n= {\rm diag}(A_n)-A_n$, then $Q_n \geq 0$ and
\begin{equation}\label{exs7.2}
 A_n^{-1} \ = \ \sum_{k=0}^\infty \: \big ( K_n^{-1} Q_n \big)^k\,
 K_n^{-1}. 
\end{equation}
Let ${\bf u}_n(r)$ be the approximation of ${\bf u}_n$  after
$r$ iterations. By taking the stopping criteria to be 
$\lnorm {\bf u}_n(r+1) - {\bf u}_n(r) \lnorm_1 < 10^{-9}$, 
we have found that the iteration terminates after $r=182000$ iterations.
Then we compared the numerical approximations ${\bf u}_n(210)$ and values of the solution 
in $l_\infty(G_n(D) \setminus \{\mst\})$ and $l_1(G_n(D) \setminus \{\mst\})$--norms, 
where the set $G_n(D) \setminus \{\mst\}$ is chosen naturally instead of the set $G_n$
because the solution $u$ of this example is not defined at $\mst$. Let the
corresponding norms be denoted by $\lnorm \cdot \lnorm_p', p =1,\infty$.
The solutions are compared according to the relative error:
\begin{equation}\label{exs7.3}
 \varepsilon_{p,rel} \ = \ \frac{\lnorm {\bf u}^\ast \,-\,{\bf u}_n(210)\lnorm_p'}{
 \lnorm {\bf u}^\ast \lnorm_p'}, \quad p = 1, \infty,
\end{equation}
where $\big({\bf u}^\ast\big)_{\mskd} = u(h\msk)$. We obtained 
$\varepsilon_{1,rel} = 0.003371$ and $\varepsilon_{\infty,rel} = 1.6254$. The value 
$\varepsilon_{\infty,rel} = 1.6254$ is realized at the grid-knot $\mst - h\mse_1$,
i.e. at one of the nearest grid-knots to the singular point $\mst$ of
solution. Actually, the difference between two solutions at this grid-knot
is $0.0075$. Let us mention that the respective error 
$\varepsilon_{1,rel} = 0.267$ is obtained with the grid-step $h=1/200$.

\begin{example}\label{exam7.2}{\rm 
Next, let us consider a differential operator $A = -\sum_{ij = 1}^2 \partial_i
a_{ij}\partial_j$ with the diffusion tensor 
\[ a \ = \ \left [ \begin{array}{cc} \sigma^2 & \alpha(\msx) \\ \alpha(\msx) & 1
 \end{array} \right ],  \quad \alpha(\msx) \:=\:\rho {\bbJ}_{D_0}(\msx), \quad
 \rho^2 < \sigma^2,\]
where $\sigma^2$ is a positive number, $\rho$ is a real number and $D_0 =
(1/4,3/4)^2$ (see Figure \ref{fgn7.1}).
The function $\msx \mapsto u^\ast(x_1,x_2) = x_1x_2$ is the unique solution 
to the boundary value problem
\[\begin{array}{lll}
 \big(A u \big)(\msx)\ = \ \mu(\msx) &{\rm for}& \msx \in D,\\
 u|\partial D \ = \  u^\ast|\partial D, &&\end{array} \]
where
\[ \mu(\msx) \ = \  2\,\rho \,{\bbJ}_{D_0}(\msx) \,+\,\frac{\rho}{4}
 \Big[ \delta \Big(x_1-\frac{1}{4}\Big)- 3\delta\Big(
 x_1-\frac{3}{4}\Big)+ \delta\Big(x_2- \frac{1}{4}\Big)-3 \delta\Big(x_2-
 \frac{3}{4}\Big)\Big].\]
The set ${\bbR}^2$ is discretized as in the previous example.
The sets $G_n(2) = \overline{D_0}\cap G_n$ and $G_n(1) = G_n \setminus 
G_n(0)$ define a partition of $G_n$. Let $G_n(1,D) = G_n(1) \cap D$.
Then $G_n(2), G_n(1,D)$ is a partition of $G_n(D)$ to be used in constructions
of numerical grids and approximate solutions.
}\end{example}

To demonstrate the efficiency of the extended schemes we choose 
$\sigma^2 = 10, \rho =2$ and the scheme parameters 
$r_1 = 1, r_2 =3$ as illustrated in Figure~\ref{fig4.2}.
These values of parameters ensure the compartmental structure of 
the system matrix. The linear system~(\ref{ex6.1}) is solved iteratively as in the previous example. 
Let ${\bf u}_n(r)$ be the approximation of ${\bf u}_n$  after
$r$ iterations. By taking the stopping criteria to be 
$\lnorm {\bf u}_n(r+1) - {\bf u}_n(r) \lnorm_1 < 10^{-9}$, 
we have found that the iteration terminates after $r=210$ iterations.
Then we compared the numerical approximations ${\bf u}_n(210)$ and 
values of the solution in $l_\infty(G_n(D))$ and $l_1(G_n(D))$--norms. 
We obtained $\varepsilon_{1,rel} = 0.8\times 10^{-5}$.
In addition, $\lnorm {\bf u}_n^\ast-{\bf u}_n(210)\lnorm_{\infty,rel} = 0.4
\times 10^{-2}$ is realized at the grid-knot with coordinates
$\msx = (0.75, 0.75)$.
\section{DISCUSSION AND CONCLUSION}\label{sec9}
Attempts to discretize a second order elliptic differential operator by 
matrices with compartmental structure go back to the work of Motzkyn and 
Wasov \cite{MW}. Their construction is based on rotations, which
can be avoided in $2$-dimensional cases by using extended schemes 
which we propose in this work. These can be used 
to get monotone schemes for any 2-dimensional problem with the second order
elliptic operator, in divergence or non-divergence form. For operators
in divergence form one has to use the construction of system matrix described 
in Section~\ref{sec4}, while for operators in non-divergence form 
system matrices are of a simpler structure. Entries of system matrices
are linear functions of $a_{ij}(\msv)$, where $\msv$ is the considered grid-knot.

In the case of dimension $d \geq 3$ the construction of Section \ref{sec4}
does not always produce monotone schemes. An additional condition 
on the matrix-valued function $\msx \mapsto \{a_{ij}(\msx)\}_{11}^{dd}$ which
is described in Assumption~\ref{ass4.1} assures the compartmental
structure of matrices $A_n$.
\begin{lemma} Let $A_n$ on $G_n$ be consistent discretizations of the differential 
operator $A(\msx) = -\sum_{ij} \partial_i a_{ij}(\msx) \partial_j$ and 
let the entries $\big(A_n\big)_{\mskd \msld}$ be linear combinations of
$a_{ij}(\msz)$ with some $\msz \in S_n(\msp(l),h\msx)$. Then $A_n$ can
possess the compartmental structure iff the tensor-valued function $\hat{a}$
is strictly positive definite on ${\bbR}^d$.
\end{lemma}
{\Proof} For a constant diffusion tensor we prove easily this lemma by
demanding the following property at each grid-knot:
\[ \big(A_n {\bf f}_n \big)_{\mskd} \ = \ A(\msx)f(\msx), \quad
 \msx = h \msk,\]
where $\msx \mapsto f(\msx)$ is any second degree polynomial in
variables $x_i$. This property is equivalent to the consistency.
To prove the necessary part, we assume that $\hat{a}(\msx_0)$ is negative definite and $A_n$
have the compartmental structure. Then the same assertion must
hold for a tensor-valued function which is constant locally at $\msx_0$. 
This contradicts the first step. {\QED}

Hence, if the matrix-valued function $\hat{a}$ is indefinite
the only way out is to use rotations. 
Unfortunately undesired features may appear;
in a subdomain $D_l \subset D$ where a rotation is required,
some larger diagonal submatrices can be reducible.
This causes the surfaces $u(n,\msx)=const$ to have
a saw-like behavior. In this case an averaging procedure of grid-solutions
helps to obtain reasonable results.

A different approach to the construction of monotone schemes for a class of 
elliptic operators was presented in Samarskii {\em et al}.~\cite{SMMM}. They constructed
schemes with various grid-steps along coordinate axis which we called basic 
schemes in our analysis. In this way they obtained monotone schemes for
a wide class of elliptic operators.

\end{document}